\documentclass{amsart}

\usepackage{amsmath,amsthm}
\usepackage{amsfonts,amssymb}

\newtheorem{thm}{Theorem}[section]
\newtheorem{cor}[thm]{Corollary}
\newtheorem{lem}[thm]{Lemma}
\newtheorem{prop}[thm]{Proposition}

\newtheorem{algo}[thm]{Algorithm}

\theoremstyle{remark}

 \def\xb{{\mathbf x}}
 \def\yb{{\mathbf y}}

 \def\CA{{\mathcal A}}     
 \def\CB{{\mathcal B}}

 \def\CH{{\mathcal H}}
 \def\CI{{\mathcal I}}

 \def\CR{{\mathcal R}}     
 \def\CS{{\mathcal S}}     
      
 \def\CV{{\mathcal V}}

 \def\NN{{\mathbb N}}
 \def\PP{{\mathbb P}}
 
 \def\RR{{\mathbb R}}

        \def\proj{\operatorname{proj}}

\newcommand{\wt}{\widetilde}
\newcommand{\wh}{\widehat}

\begin{document}
\input epsf

\title[Attenuated Radon projections]
{Approximation and Reconstruction from Attenuated Radon Projections} 
 
\author{Yuan Xu}
\address{ Department of Mathematics\\ University of Oregon\\
    Eugene, Oregon 97403-1222.}
\email{yuan@math.uoregon.edu}
\author{Oleg Tischenko}
\author{Christoph Hoeschen}
\address{Institute of Radiation Protection\\
GSF - National Research Center for Environment and Health\\
D-85764 Neuherberg, Germany}
\email{oleg.tischenko@gsf.de}
\email{christoph.hoeschen@gsf.de}

\date{\today}

\begin{abstract}
Attenuated Radon projections with respect to the weight function
$W_\mu(x,y) = (1-x^2-y^2)^{\mu-1/2}$ are shown to be closely related to 
the orthogonal expansion in two variables with respect to $W_\mu$. 
This leads to an algorithm for reconstructing two dimensional functions 
(images) from attenuated Radon projections. Similar results 
are established for reconstructing functions on the sphere from projections 
described by integrals over circles on the sphere, and for reconstructing 
functions on the three-dimensional ball and cylinder domains. 
\end{abstract}

\maketitle

\section{Introduction}
\setcounter{equation}{0}
 
Computer tomography (CT) offers a non-invasive method for 2D cross-sectional
or 3D imaging of an object. In a typical CT application, the distribution of 
the attenuation coefficient through a body from measurements of x-ray 
transmission is estimated and used to reconstruct an image of the object. 
The mathematical foundation of CT is Radon transform. Let $f$ be a function 
defined on the unit disk $B^2$ of the $\RR^2$ plane. A Radon transform of $f$ 
is a line integral, 
\begin{equation} \label{eq:0.1}
\CR_\theta (f;t) :=  \int_{I(\theta,t)} f(x,y) dx dy,
  \qquad 0 \le \theta \le 2\pi,   \quad -1\le t \le 1, 
\end{equation}
where $I(\theta,t) = \{(x,y): x\cos\theta+y\sin \theta 
= t\} \cap B^2$ is a line segment inside $B^2$. An essential problem in CT 
is to reconstruct the function $f$ from its Radon projections. An algorithm 
amounts to an approximation to $f$ that uses values of $\CR_\theta(f;t)$ from 
a finite set of parameters $(\theta, t)$. 

The attenuation of an x-ray beam is dependent on the energy of each photon.
A line integral as defined in \eqref{eq:0.1} represents a monochromatic x-ray.
In practice, however, an x-ray is usually polychromatic, meaning that it 
consists of photons with different energies. This could lead to artifacts in 
the reconstruction; see, for example, \cite[Chapt. 4]{KS}. A polychromatic 
x-ray is represented by the so-called attenuated Radon projections for which 
the integral is taken against $\exp\{-\alpha_\theta(x,y)\} dxdy$, where 
$\alpha_\theta(x,y)$ is a given function, instead of $dxdy$. Attenuated Radon 
transform appears in, for example, emission tomography \cite{N}. The 
reconstruction algorithms for attenuated Radon data have been derived from
Novikov's inversion formula (\cite{No} and \cite{N2}). See also the recent
survey in \cite{F} in this direction.  

In the present paper we consider the special case that 
$\exp\{-\alpha_\theta(x,y)\}$ is given, or can be approximated, by the 
function 
\begin{equation} \label{eq:W-mu}
  W_\mu(x,y) = (1-x^2-y^2)^{\mu -1/2}, \qquad (x,y) \in B^2, 
\end{equation}
where $\mu \ge 0$; in other words, $\alpha_\theta(x,y) = - (\mu-1/2)
 \log (1-x^2-y^2)$. The attenuated Radon transform, denote by 
$\CR_\theta^\mu$, then takes the form 
\begin{equation} \label{eq:R-mu}
\CR_\theta^\mu(f;t) :=  \int_{I(\theta,t)} f(x,y) W_\mu(x,y) dx dy,
  \qquad 0 \le \theta \le 2\pi,   \quad -1\le t \le 1. 
\end{equation}
Clearly this is just a special case of the attenuated Radon transform.
This case, however, appears to be useful in understanding the effect of 
monochromatic and polychromatic x-rays. In this regard let us mention 
the classical example of the water phantom in a skull in \cite[p. 121]{KS},
which demonstrated that beam hardening causes an elevation in CT numbers
for tissues close to the skull bone. The attenuated Radon transform 
defined in \eqref{eq:R-mu} models the boundary behavior of the x-rays
differently. 

Our approach is based on orthogonal polynomial expansions on $B^2$. Let 
$\CV_n^2(W_\mu)$ denote the space of orthogonal polynomials with respect to 
the weight function $W_\mu$ on $B^2$. It is well known that 
$$
  L^2(B^2,W_\mu) = \sum_{k=0}^\infty \bigoplus \CV_k^2: \qquad 
       f = \sum_{k=1}^\infty \proj_k^\mu f,  
$$
where $\proj_k^\mu f$ is the projection of $f$ on $\CV_k^2(W_\mu)$. The 
infinite series holds in the sense that the sequence of the partial sums 
$$
   S_n^\mu(f; x,y): = \sum_{k=0}^n \proj_k^\mu f(x,y), \qquad n \ge 0,  
$$
converges to $f$ as $n \to \infty$ in $L^2(B^2,W_\mu)$ norm. The partial sum 
$S_n f$ provides a natural approximation to $f$. It turns out that there is 
a remarkable connection between $S_n^\mu f$ and the attenuated Radon 
transforms, which states that 
\begin{equation} \label{eq:partial}
S_{2m}f(x,y) = \sum_{\nu = 0}^{2m} \int_{-1}^1\CR_{\phi_\nu}^\mu(f;t)
 \Phi_\nu(t;x,y) dt,  \qquad \phi_\nu = \frac{2\nu \pi}{2m+1},
\end{equation} 
where $\Phi_\nu$ are polynomials of two variables given by explicit formulas.
This representation provides a simple and direct access to attenuated Radon 
data. For the ordinary Radon transforms ($\mu = 1/2$), this was 
discovered recently in \cite{X05}. Applying an appropriate quadrature formula 
to the integrals in the expression leads to an approximation to $f$ that uses 
discrete attenuated Radon projections. One important feature of the algorithm
is that polynomials up to a certain degree are reconstructed exactly, which
guarantees that the algorithm has a fast rate of convergence. Such an algorithm
can be easily implemented numerically. For the ordinary Radon transforms, the
algorithm is named OPED (Orthogonal Polynomial Expansion on the Disk) and it
has proved to be a highly effective method \cite{XTC,XTC2}. 

There are other expressions in the spirit of \eqref{eq:partial}. In order to
prove them, we need to study orthogonal expansions in terms of orthogonal 
polynomials with respect to $W_\mu(x,y)$ on $B^2$. The case $\mu =1/2$ is
easier since an orthonormal basis for $\CV_k^2(W_{1/2})$ is known to be 
$U_k(x \cos \tfrac{j \pi}{k+1}+ y \sin \tfrac{j \pi}{k+1})$, $0 \le j \le k$. 
No such convenient orthonormal basis is available for $\mu \ne 1/2$. 

There is another advantage for considering the attenuated Radon transform
$\CR^\mu_\theta (f;t)$. It is known that there is a close relation 
between orthogonal polynomials on the unit ball and those on the unit sphere, 
which allows us to establish analogous results on the unit sphere $S^2$. In 
particular, the case $\mu =0$ on $B^2$ can be used to show that we can 
reconstruct a function $f$ from its integral projections 
\begin{equation} \label{eq:1.4}
 Q f(\zeta;t) =  \int_{\langle \xb,\zeta\rangle = t} f(\xb) d \omega(\xb),  
 \qquad 0 \ne \zeta \in S^2, \quad  -1\le t \le 1, 
\end{equation}
where $\xb = (x_1,x_2,x_3)$ and $d\omega$ is the surface measure on $S^2$. 
Reconstruction from such spherical transforms has been studied in the 
literature, see \cite{NW}. 

From the disk $B^2$ we can also extend the results to the unit 
ball $B^3$ and to cylinder domains in $\RR^3$, taking Radon projections on
parallel disks in each case. It turns out, however, that there is an important
difference between the ball and the cylinder. For the cylinder domain, all 
results obtained in the disk can be extended without problem. For the unit
ball, however, we still have an analogue of \eqref{eq:partial} but the 
reconstruction algorithm may no longer work as efficient as in the cylinder 
case. The problem is that the operator produced by the algorithm no longer 
preserves polynomials. 

For the algorithm on $B^2$, we provide a numerical example in Section 2, 
which reconstructs a 2D phantom image for three different values of $\mu$. 
For the transform on the sphere and the 3D transforms on the ball and on
the cylinder domain, we will content with deriving the algorithms and will 
not discuss convergence or the performance of the algorithms at this time. 
 
The paper is organized as follows. In the following section we consider the
reconstruction and approximation on the unit disk $B^2$ from attenuated
Radon projections. This section is divided into several subsections, the
last one includes the numerical example. In Section 3 the results on $B^2$ 
are transplanted to those on the surface $S^2$, 
while the attenuated Radon projections become weighted spherical transforms.
The analogous results are then established for the unit ball $B^3$ in 
Section 4 and for the cylinder domain in Section 5. 
 
\section{Reconstruction and Approximation on the unit disk}
\setcounter{equation}{0}

Let $\Pi^d$ denote the space of polynomials of $d$ variables and let
$\Pi_n^d$ denote the subspace of polynomials of total degree $n$ in 
$\Pi^d$, which has dimension $\dim \Pi_n^d = \binom{n+d}{d}$. We set
$\Pi_n:=\Pi_n^1$. In 
this section we mainly work with the case $d =2$. 

\subsection{Orthogonal polynomials on the unit disk}

Let $W_\mu$ be the weight function defined in \eqref{eq:W-mu}.
Let $\CV_k^2(W_\mu)$ denote the space of orthogonal polynomials of 
degree $k$ on $B^2$ with respect to the inner product 
$$
   \langle P, Q \rangle_\mu = a_\mu \int_{B^2} P(x,y) Q(x,y)W_\mu(x,y) dx dy, 
      \qquad a_\mu = (\mu+1/2)/\pi, 
$$
where $a_\mu$ is the normalization constant of $W_\mu$, $a_\mu = 
1/\int_{B^2}W_\mu(x) dx$. Thus, $P \in \CV_k^2(W_\mu)$ if $P$ is of degree 
$k$ and $\langle P, Q\rangle_\mu = 0$ for all $Q \in \Pi_{k-1}^2$. We note that
elements in a basis for $\CV_k^2(W_\mu)$ may not be orthogonal with respect
to each other according to our definition. A basis for $\CV_k^2(W_\mu)$ is 
called orthonormal if the elements in the basis are mutually orthogonal and 
$\langle P, P\rangle_\mu =  1$. 

The reproducing kernel of the space $\CV_k^2(W_\mu)$ plays an important
role in our development. In terms of an orthonormal basis  $\{P_j^k: 0 \le j \le k\}$ 
of $\CV_k^2(W_\mu)$, the reproducing kernel satisfies
\begin{equation} \label{eq:Pkernel}
       P_k(W_\mu;\xb,\yb) = \sum_{j=0}^k P_j^k(\xb) P_j^k(\yb).  
\end{equation}
The kernel is independent of the choice of the bases of $\CV_k^2(W_\mu)$. 
In fact, a compact formula for this kernel can be given in terms of the 
Gegenbauer polynomial \cite{X99},
\begin{align} \label{eq:reprod}
 & P_k(W_\mu;\xb,\yb) =  \frac{k+\mu+1/2}{\mu+1/2} b_{\mu-1} \\
  & \qquad \times \int_{-1}^1 C_k^{\mu+1/2} \left(\langle \xb, \yb \rangle 
         +\sqrt{1-\|\xb\|^2}
        \sqrt{1-\|\yb\|^2}\,t \right)(1-t^2)^{\mu-1} dt  \notag
\end{align} 
for $\mu > 0$, the formula also holds for $\mu =0$ upon taking limit
$\mu \to 0$. Here and in the following, the Gegenbauer polynomials 
$C_k^\lambda(s)$ are orthogonal with respect to $(1-s^2)^{\lambda -1/2}$ 
on $[-1,1]$, 
\begin{equation} \label{eq:Gegenbauer}
   c_{\lambda-1/2} \int_{-1}^1 C_k^\lambda(s) C_l^\lambda(s) 
     (1-s^2)^{\lambda-1/2}ds = \frac{\lambda(2\lambda)_k}{(k+\lambda) k!}
          \delta_{k,l} : = h_k \delta_{k,l}, 
\end{equation}
where $c_{\lambda-1/2} := \Gamma(\lambda+1)/ (\sqrt{\pi} \Gamma(\lambda+1/2))$
is the normalization constant of the weight function $(1-s^2)^{\lambda-1/2}$ 
on $[-1,1]$, and $(a)_k : = a(a+1)\cdots (a+k-1)$.  For $\mu = 1/2$, 
$C_k^{\mu+1/2}(s)= U_k(s)$ is the Chebyshev polynomial of the second kind. 

For the weight function $W_{1/2}(x) = 1$, it is known \cite{LS} that the set 
$$
\left \{ U_k \left(x  \cos \theta_{j,k}  + y \sin \theta_{j,k}  \right): 
      0 \le j \le k  \right\}
$$
forms an orthonormal basis of $\CV_k^2(W_{1/2})$. The elements of this
basis are the so-called ridge functions. In general, given an angle $\phi$ 
and a polynomial $p \in \Pi_k:=\Pi_k^1$, a ridge polynomial is defined by 
$$
   p(\phi; x,y) : = p(x \cos \phi + y \sin \phi), \qquad \phi \in [0,2\pi].
$$
It is easy to see that $p(\phi;x,y)$ is a polynomial in $\Pi_k^2$ as 
well. The functions $\{C_k^{\mu+1/2}(\theta_{j,k};x,y): 0 \le j \le k\}$, where
$\theta_{j,k} = j\pi/(k+1)$, form a basis for $\CV_k^2(W_\mu)$, abeit not 
an mutually orthogonal one (see, for example, \cite{X00}). The lack of 
orthonormal ridge basis in the case of $\mu \ne 1/2$ makes the results
for attenuated Radon transform more difficult, as we shall see below.

We call a polynomial $P \in \Pi_k$ of one variable {\it symmetric} with 
respect to the origin if $P$ is even when $k$ is even, and $P$ is odd 
when $k$ is odd. It is known that $C_k^{\mu+1/2}(t)$ is symmetric with
respect to the origin. The ridge polynomials arising from such a polynomial 
turn out to satisfy a remarkable relation. 

%%%%%%%%%%%%%%%%%
%%%%  Proposition 2.2 %%%%
%%%%%%%%%%%%%%%%%

\begin{prop} \label{prop:key}
For $n \ge 0$ and $k \le n$, the identity
\begin{equation} \label{eq:sumUAll}
 \frac{1}{n+1} \sum_{\nu=0}^{n}
   U_k\left(\tfrac{\nu \pi}{n+1};\cos\theta,\sin\theta\right)
        P_k\left( \tfrac{\nu \pi}{n+1};x,y\right) = P_k(\theta;x,y)
\end{equation}
holds for all polynomials $P_k \in \Pi_k$ that are symmetric with 
respect to the origin.  
\end{prop}

\begin{proof}
The proof uses the following elementary trigonometric identities
\begin{equation}\label{eq:TrigIdent}
 \sum_{\nu=0}^{n} \sin k \tfrac{2 \nu \pi}{n+1} = 0 \quad \hbox{and} \quad 
 \sum_{\nu=0}^{n} \cos k \tfrac{2 \nu\pi}{n+1} = \begin{cases} n+1, &
     \hbox{if $k = 0 \mod n+1$} \\
        0, &\hbox{otherwise} \\
   \end{cases}
\end{equation}
that hold for all nonnegative integers $k$. Let us prove the case $k = 2l$. 
%It suffices to prove \eqref{eq:sumUAll} for $P_k(t)$ being the Chebyshev
%polynomial $U_k(t)$. 
We follow the proof of Proposition 2.3 in \cite{X05}. The polynomial $P_k$ 
can be written as a linear combination of $U_{k-2j}$ for $0 \le 2j \le k$.
Consequently, we can write $P_{2l} (\theta;x,y)$ as
\begin{align} \label{eq:u1}
  P_{2l}(\theta; x,y) = P_{2l}(r \cos(\theta - \phi))   
          =  \sum_{j=0}^l b_j(r) \cos 2j (\theta - \phi)
\end{align}
in polar coordinates $x = r \cos \phi$ and $y = r \sin \phi$, where 
$b_j(r)$ is a polynomial of degree $2 j$ in $r$. Furthermore, we know
that 
$$
 U_{2l}(\theta; \cos \phi,\sin \phi) = U_{2l}( \cos(\theta - \phi))   
     = \sum_{j=0}^l d_j \cos 2j (\theta - \phi)
$$
where $d_0 = 1$ and $d_j = 2$ for $j \ge 1$. The identities 
\eqref{eq:TrigIdent} and the product formula of the cosine function
shows that 
$$
\frac{1}{n+1}\sum_{\nu=0}^{n} \cos 2i (\theta - \tfrac{\nu \pi}{n+1}) 
    \cos 2j (\phi - \tfrac{\nu\pi}{n+1}) 
  = \begin{cases} 0, & \hbox{if $i\ne j$}, \\ 
        \frac{1}{2} \cos 2 j(\theta - \phi), & \hbox{if $0<i = j \le n$},\\
         1, & \hbox{if $i = j = 0.$} \end{cases}
$$  
Let us denote by $I_k$ the left hand side of \eqref{eq:sumUAll}. The above
trigonometric identity implies immediately that, for $0 \le 2l \le n$, 
\begin{align*}
 I_{2l} & = \sum_{i= 0}^l d_i \sum_{j =0}^l b_j(r) 
  \frac{1}{n+1} \sum_{\nu=0}^{n} \cos 2i (\theta - \tfrac{\nu \pi}{n+1})
       \cos 2j (\phi - \tfrac{\nu \pi}{n+1}) \\
   & =  \sum_{j=0}^l b_j(r) \cos 2 j(\theta - \phi) =
     P_{2l}(r \cos (\theta - \phi)) = P_{2l}(\theta; x,y).
\end{align*}
This completes the proof for the case $k = 2l \le n$.  The case $k = 2l-1$ 
is similar. 
\end{proof}

In \eqref{eq:sumUAll} the summation is over angles, $\nu \pi/(n+1)$, 
that are equally spaced in the interval $[0,\pi)$. In the case that $n$ 
is even,  the angles can be arranged as equally spaced angles 
in $[0, 2\pi]$ by using the fact that 
\begin{equation} \label{eq:CosSin}
  \cos \tfrac{(2k+1)\pi}{2m+1} = - \cos \tfrac{(2m+2k)\pi}{2m+1}  
  \quad\hbox{and}\quad
   \sin \tfrac{(2k+1)\pi}{2m+1} = - \sin \tfrac{(2m+2k)\pi}{2m+1}.  
\end{equation}
The result is the following proposition proved in \cite{X05} for
$P_k$ being the Chebyshev polynomial of the second kind.

\begin{prop} \label{prop:keyEven}
For $m \ge 0$ and $k \le 2m$, the identity
\begin{equation} \label{eq:sumU} 
\frac{1}{2m+1} \sum_{\nu=0}^{2m}
   U_k\left(\tfrac{2 \nu \pi}{2m+1};\cos\theta,\sin\theta\right)
        P_k\left( \tfrac{2 \nu \pi}{2m+1};x,y\right) = P_k(\theta;x,y)
\end{equation}
holds for all polynomials $P_k \in \Pi_k$ that are symmetric with 
respect to the origin.  
\end{prop}

There are many orthonormal bases of $\CV_k^2(W_\mu)$ that are known explicitly 
(see \cite{DX}). One that is particularly useful for us is given in terms 
the polar coordinates 
$$
 x = r \cos \phi, y = r \sin \phi, \qquad 0 \le r \le 1, 
\quad 0 \le \phi \le 2 \pi,
$$
and Jacobi polynomials \cite[Prop. 2.3.1]{DX}. Let $p_n^{(\alpha,\beta)}(t)$ 
denote the orthonormal Jacobi polynomials, that is, 
$$
  c_{\alpha,\beta}\int_{-1}^1 p_n^{(\alpha,\beta)}(t)p_m^{(\alpha,\beta)}(t) 
    (1-t)^\alpha(1+t)^\beta  dt = \delta_{m,n}, \quad 
          m, n = 0,1,2, ... 
$$
where $c_{\alpha,\beta}$ is the normalized constant so that
$c_{\alpha,\beta} \int_{-1}^1 (1-t)^\alpha (1+t)^\beta dt =1$.  

\begin{prop} \label{lem:2.3}
For $\varepsilon = 0$ or $1$, define the polynomials $P_{l,\varepsilon}^k$ 
by
\begin{equation}\label{eq:on-basis}
P_{l,\varepsilon}^k(x,y) = h_{l,k} p_l^{(\mu-\frac{1}{2},k-2l)}(2r^2 -1) 
    r^{k-2l} S_{k-2l,\varepsilon}(\phi),
\end{equation}  
where 
\begin{align*}
& S_{k-2l, 0}(\phi) = \cos (k-2l)\phi \quad \hbox{ for } \quad 0 \le 2l \le k,
 \\
& S_{k-2l, 1}(\phi) = \sin \,(k-2l)\phi \quad \hbox{for} \quad 0 \le 2l \le 
k-1,
\end{align*}
%$S_{k-2l, 0}(\phi) = \cos (k-2l)\phi$ for $0 \le 2l \le k$, 
%$S_{k-2l, 1}(\phi) = \sin (k-2l)\phi$ for $0 \le 2l \le k-1$, 
and 
$$
[h_{l,k}]^2  :=  \frac{ \Gamma(k-2l+\mu+3/2)}{\Gamma(\mu+3/2) \Gamma(k-2l+1)}.
$$
Then these polynomials form an orthonormal basis for $\CV_k^2 (W_\mu)$.
\end{prop}

By the definition of the reproducing kernel \eqref{eq:Pkernel} and the formula
\eqref{eq:reprod}, it follows that the above orthonormal basis satisfies
\begin{equation}\label{eq:OPsum}
 \sum_{\varepsilon=0,1}
   \sum_{0 \le 2l \le k} P_{l,\varepsilon}^k(x,y)P_{l,\varepsilon}^k(\cos \phi,\sin\phi)
    = \frac{k+\lambda}{\lambda}C_k^\lambda(\phi;x,y),
\end{equation}
where $\lambda = \mu + 1/2$. This formula will play an important role below. 
It shows, in particular, that the expansion of $C_k^{\mu+1/2}(\phi;x,y)$ in terms
of our orthonormal basis.  The following lemma shows the converse. 

\begin{lem} \label{lem:sum1}
Let $\theta_{j,k} = j \pi/(k+1)$. Then for $0 \le 2l \le k$ if $\varepsilon =0$
and $0 \le 2l \le k-1$ if $\varepsilon =1$,
$$
 \frac{1}{k+1} \sum_{j=0}^k S_{k-2 l, \varepsilon}(\theta_{j,k}) 
     C_k^{\mu+1/2} (\theta_{j,k};x,y) 
        = \frac{\mu+\frac{1}{2}}{k+\mu+\frac{1}{2}}
               H_{l,k}^\mu d_{l,k} P_{l,\varepsilon}^k(x,y),
$$
where $d_{l.k} = 1/2$ if $2l<k$ and $d_{l,k} = 1$ if $2l = k$, 
$H_{l,k}^\mu:=h_{l,k}^\mu p_l^{(\mu+1/2,k-2l)}(1)$ and 
$$
   \left[H_{l,k}^\mu \right]^2 = \frac{ (\mu+\frac12)_l 
       (\mu+ \frac{3}{2})_{k-l} (k+\mu+\frac{3}{2})}
         { l! (k-l)! (k-l+\mu+\frac{3}{2})}.
$$
\end{lem}

\begin{proof}
Using the identities \eqref{eq:TrigIdent} it is easy to verify that 
\begin{equation} \label{eq:sin-cos}
   \frac{1}{k+1} \sum_{j=0}^k S_{k-2l,\varepsilon} (\theta_{j,k} )
     S_{k-2 l',\varepsilon} (\theta_{j,k} ) =  d_{l,k} \delta_{l,l'}.
\end{equation}
Using \eqref{eq:on-basis} and the fact that 
$P_{l,\varepsilon}^k (\cos \theta_{j,k},\sin \theta_{l,k})
   = H_{l,k}^\mu S_{k-2l,\varepsilon}(\theta_{j,k})$, we obtain
\begin{align*}
&  \frac{1}{k+1} \sum_{j=0}^k S_{k-2 l, \varepsilon}(\theta_{j,k}) 
     C_k^{\mu+1/2} (\theta_{j,k};x,y) \\
 &   = \frac{\mu+\frac{1}{2}}{k+\mu+\frac{1}{2}}
   \sum_{0\le l \le 2k} P_{l,\varepsilon}^k(x,y)
   \frac{1}{k+1} \sum_{l=0}^k P_{l,\varepsilon}^k(\cos\theta_{j,k},\sin \theta_{j,k})
                S_{k-2 l, \varepsilon}(\theta_{j,k}) \\
& = \frac{\mu+\frac{1}{2}}{k+\mu+\frac{1}{2}} H_{l,k}^\mu d_{l,k} P_l^k(x,y)
\end{align*}
upon using the equation \eqref{eq:sin-cos}. Finally, the expression of
$[H_{l,k}^\mu]^2$ is derived from the well-known formula of 
$p_l^{\alpha,\beta}(1)$ (see \cite{Sz}) and the formula of $h_{l,k}^\mu$.
\end{proof}

\begin{lem} \label{lem:sum2}
Let $\theta_{j,k}$ be as above. Then 
$$
 \frac{1}{k+1} \sum_{j=0}^k S_{k-2l,\varepsilon}(\theta_{j,k})
   U_k(\theta_{j,k}; \cos \phi,\sin\phi) = d_{l,k} S_{k-2l,\varepsilon}(\phi).  
$$
\end{lem} 

\begin{proof}
Using \eqref{eq:u1} and the analog formula for $U_{2l-1}$, the identity
is an easy consequence of \eqref{eq:sin-cos}.
\end{proof}

\subsection{Attenuated Radon transforms}

Let $\theta$ be an angle  measured counterclockwise from the positive 
$x$-axis. Let $\ell$ denote the line perpendicular to the direction
$(\cos \theta,\sin\theta)$ and passes through the point $(t\cos
\theta, t\sin \theta)$. The equation of the line is $\ell(\theta,t) = \{(x,y):
 x \cos \theta + y \sin \theta = t\}$ for $-1 \le t \le 1$.  We use
 \begin{equation} \label{eq:line}
I(\theta, t) = \ell (\theta, t) \cap B^2, \qquad 
      0 \le \theta< 2 \pi, \quad -1 \le t \le 1,
\end{equation}
to denote the line segment of $\ell$ inside $B^2$.  Let $W_\mu$ be the 
weight function defined in \eqref{eq:W-mu}. The attenuated Radon 
projection of a function $f$, with respect to $W_\mu$, in the 
direction $\theta$ with parameter $t \in [-1,1]$ is defined in 
\eqref{eq:R-mu}. It can be written as 
\begin{align} \label{eq:Radon}
\CR_\theta^\mu(f;t) = \int_{-\sqrt{1-t^2}}^{\sqrt{1-t^2}}  
    f(t \cos \theta - s \sin \theta, t\sin\theta+ s \cos \theta) 
       W_\mu(s,t) ds,
\end{align}
using the fact that the mapping $(s,t) \mapsto (x,y)$ defined by 
$x = t \cos \theta - s \sin \theta$ and 
$y= t\sin\theta+ s \cos \theta$ amounts to a rotation. When $\mu =1/2$,
this is the usual Radon projection, which is also called an X-ray transform.
The definition \eqref{eq:R-mu} or \eqref{eq:Radon} shows that 
$\CR_\theta^\mu(f;t)= \CR_{\pi+\theta}^\mu(f;-t)$. 
 
The ridge polynomials are particularly useful for studying Radon transforms, 
as seen in the following result:

%%%%%%%%%%%%%%%%%
%%%%  Proposition 2.1 %%%%
%%%%%%%%%%%%%%%%%

\begin{prop} \label{prop:2.1}
For $f \in L^1(B^2)$ and $p \in \Pi_k$, 
\begin{equation} \label{eq:2.2}
 \int_{B^2} f(x,y) p(\phi; x,y) W_\mu(x,y)dx dy = 
    \int_{-1}^1 \CR_\phi^\mu (f;t)  p(t) dt. 
\end{equation}
\end{prop}

\begin{proof}
Since the change of variables $t = x \cos \phi + y \sin \phi$ and 
$s =-x \sin \phi + y \cos\phi$ amounts to a rotation, we have 
\begin{align*}
 & \int_{B^2} f(x,y) p_k(\phi; x,y)W_\mu(x,y) dx dy  \\
 & \qquad = 
    \int_{B^2} f(t \cos \phi - s \sin \phi, t \sin \phi + s \cos \phi)  
     p_k(t) W_\mu(t,s) dt ds \\
 & \qquad  = \int_{-1}^1 \int_{-\sqrt{1-t^2}}^{\sqrt{1-t^2}}
  f(t\cos\phi-s\sin\phi,t\sin \phi + s \cos \phi) W_\mu(t,s)ds p_k(t) dt, 
\end{align*}  
the inner integral is exactly $\CR_\phi^\mu(f;t)$ by \eqref{eq:Radon}. 
\end{proof}

In particular, attenuated Radon transforms of the orthogonal polynomials 
in $\CV_n^2(W_\mu)$ can be explicitly computed. 

%%%%%%%%%%%%%%%%%
%%%%%%  Lemma 2.4 %%%%
%%%%%%%%%%%%%%%%%

\begin{lem} \label{lem:2.4}
If $P \in \CV_k^2(W_\mu)$ then for each $t\in (-1,1)$, $0\le\theta \le 2\pi$, 
\begin{equation} \label{eq:Marr}
 \CR_\theta^\mu (P ; t) =  b_\mu (1-t^2)^\mu 
    \frac{C_k^{\mu+1/2}(t)}{C_k^{\mu+1/2}(1)}P(\cos \theta, \sin \theta),
\end{equation}
where $b_\mu = c_\mu^{-1}$ for $c_\mu$ defined in \eqref{eq:Gegenbauer}.
\end{lem} 

\begin{proof}
Changing variables in \eqref{eq:Radon} shows that 
\begin{align*}
Q(t): = &  (1-t^2)^{-\mu}\CR_\theta^\mu (P; t ) \\
  =  & \int_{-1}^1 
    P \left(t \cos \theta - s \sqrt{1-t^2} \sin \theta, t \sin \theta + 
        s \sqrt{1-t^2} \cos \theta\right) (1-s^2)^{\mu-1/2} ds.
\end{align*}
Since an odd power of $\sqrt{1-t}$ in the integrand is always attached 
with an odd power of $s$,  which has integral zero, $Q(t)$ is a polynomial 
of $t$ of degree at most $k$.  Furthermore, the integral shows that $Q(1) =
b_\mu P(\cos \theta, \sin\theta)$. The equation \eqref{eq:2.2} in 
Proposition \ref{prop:2.1} shows that
$$
 \int_{-1}^1  \frac{\CR_\theta^\mu (P; t )} {(1-t^2)^{\mu}} C_j^{\mu+1/2}(t) 
      (1-t^2)^\mu dt  =  \int_{B^2} P(x,y) C_j^{\mu+1/2}(\theta; x,y) dxdy = 0,
$$
for $ j=0,1,\ldots, k-1$, since $P \in \CV_k(B^2)$. In particular, this 
shows that $Q(t)$ is in fact orthogonal to all polynomials in $\Pi_{k-1}$
with respect to the weight function $(1-t^2)^\mu$ on $[-1,1]$. Since $Q$ is 
of degree $k$, it must be an orthogonal polynomial of degree $k$ with 
respect to this weight function. Hence, we conclude that  $Q(t)= c 
C_k^{\mu+1/2}(t)$ for some 
constant $c$ independent of $t$. Setting $t =1$ shows that $c = 
b_\mu  P(\cos \theta,\sin\theta)/C_k^{\mu+1/2}(1)$. 
\end{proof}

In the case of $\mu = 1/2$, the above lemma appeared first in \cite{Marr}. 

\subsection{Orthogonal expansion and attenuated Radon
projections}

The standard Hilbert space theory shows that any function in $L^2(W_\mu;B^2)$
can be expanded as a Fourier orthogonal series in terms of $\CV_n^2(W_\mu)$. 
More precisely,  
\begin{equation} \label{eq:fourier}
  L^2(W_\mu;B^2) = \sum_{k=1}^\infty \bigoplus \CV_k^2(W_\mu): 
     \qquad   f = \sum_{k=1}^\infty   \proj_k^\mu f, 
\end{equation} 
where $\proj_k^\mu f$ is the orthogonal projection of $f$ from 
$L^2(W_\mu;B^2)$ onto the subspace $\CV_k^2(W_\mu)$.  It is well known 
that $\proj_k^\mu f$ can be written as 
an integral operator in terms of the reproducing kernel 
$P_k(W_\mu;\cdot,\cdot)$ of $\CV_k(B^2)$ in $L^2(B^2)$; that is, 
\begin{equation} \label{eq:Pintegral}
   \proj_k^\mu f(\xb) = \int_{B^2} P_k(W_\mu;\xb,\yb) f(\yb)  W_\mu(\yb)d\yb, 
\end{equation}
where $\xb = (x_1,x_2)$ and $\yb = (y_1,y_2)$. 

This formula plays an essential role in studying the convergence 
behavior of the orthogonal expansions, see for example \cite{X99,X01}.  
For our purpose, we need a different expression for $\proj_k f$. This 
is the following remarkable formula that relates $\proj_k f$ to the 
attenuated Radon transforms of $f$ directly.  Let 
$$
     \xi_\nu = \frac{\nu \pi }{n+1}, \qquad 0 \le \nu \le n.
$$
    
%%%%%%%%%%%%%%%%%
%%%%%%  Theorem 2.5 %%%
%%%%%%%%%%%%%%%%%

\begin{thm}  \label{thm:proj-sum}
For $n \ge 0$ and $k \le n$,  the operator $\proj_k^\mu f$ can be 
written as
\begin{align} %\label{eq:proj} 
\proj_{k}^\mu f(x,y) & = \frac{1}{n+1}  \sum_{\nu=0}^{n} 
  a_\mu \int_{-1}^1 \CR_{\xi_\nu}^\mu(f;t) D_k^{\mu+1/2}(\xi_\nu,t; x,y) dt
  \label{eq:proj}  \\ 
& = \frac{1}{2n+2}  \sum_{\nu=0}^{2n+1} 
  a_\mu \int_{-1}^1 \CR_{\xi_\nu}^\mu(f;t) D_k^{\mu+1/2}(\xi_\nu,t; x,y) dt
\label{eq:proj2}   
\end{align}  
where 
\begin{equation} \label{eq:C-t-xi} 
     D_k^{\mu+1/2}(\xi, t; x,y) =   \frac{k+\mu+1/2}{\mu+1/2} 
       C_k^{\mu+1/2}(t) D_k^{\mu+1/2}(\xi;x,y)
\end{equation}
with $\lambda_{l,k}^\mu = [H_{l,k}^\mu]^{-2}$ and 
$$ 
 D_k^{\mu+1/2}(\xi_\nu;x,y) :=
    \sum_{l=0}^k  \lambda_{l,k}^\mu P_l^k(\cos \xi_\nu,\sin\xi_\nu)P_l^k(x,y). 
 $$
\end{thm}

\begin{proof}
Since $C_k^{\mu+1/2}$ is symmetric with respect to the origin, using 
Proposition \ref{prop:key} and Proposition \ref{prop:2.1}, we have 
\begin{align*}
 &   a_\mu \int_{B^2} f(x,y) C_k^{\mu+1/2}(\theta_{j,k};x,y) W_\mu(x,y)dxdy \\
 & \qquad =  \frac{1}{n+1} \sum_{\nu=0}^{n} 
    U_k (\xi_\nu; \cos \theta_{j,k},\sin\theta_{j,k}) \\
  & \qquad \qquad  \times  a_\mu \int_{B^2} f(x,y)
    C_k^{\mu+1/2}(\xi_\nu ;x,y) W_\mu(x,y)dxdy \\
  & \qquad =  \frac{1}{n+1} \sum_{\nu=0}^{n} 
         a_\mu \int_{-1}^1  \CR_{\xi_\nu}(f;t) C_k^{\mu+1/2}(t) dt 
          U_k(\xi_\nu; \cos \theta_{j,k}, \sin \theta_{j,k}). 
\end{align*}
Using Lemma \ref{lem:sum1} and Lemma \ref{lem:sum2} we conclude that 
\begin{align*}
 &   a_\mu \int_{B^2} f(x,y) P_{l,\varepsilon}^k(x,y) W_\mu(x,y)dxdy \\
  & \qquad =  \frac{1}{n+1} \sum_{\nu=0}^{n} 
         a_\mu \int_{-1}^1  \CR_{\xi_\nu}(f;t) C_k^{\mu+1/2}(t) dt 
           \frac{k+\mu+\frac{1}{2}}{\mu+\frac{1}{2}}  [ H_{l,k}^\mu]^{-1} \\
    & \qquad\qquad \times d_{l,k}^{-1} \frac{1}{k+1} \sum_{j=0}^k 
           S_{k-2l,\varepsilon}(\theta_{j,k})
            U_k(\xi_\nu; \cos \theta_{j,k}, \sin \theta_{j,k})\\        
    & \qquad =  \frac{1}{n+1} \sum_{\nu=0}^{n} 
         a_\mu \int_{-1}^1  \CR_{\xi_\nu}(f;t) C_k^{\mu+1/2}(t) dt 
           \frac{k+\mu+\frac{1}{2}}{\mu+\frac{1}{2}}  [ H_{l,k}^\mu]^{-1} S_{k-2l,\varepsilon} 
              (\xi_\nu).         
\end{align*}
Multiplying by $P_{l,\varepsilon}^k (x,y)$ and sum up, it follows from 
the definition of the reproducing kernel that 
\begin{align*}
\proj_k^\mu f(x,y) = & \frac{1}{n+1} \sum_{\nu=0}^{n} 
         a_\mu \int_{-1}^1  \CR_{\xi_\nu}(f;t) C_k^{\mu+1/2}(t) dt 
 \frac{k+\mu+\frac{1}{2}}{\mu+\frac{1}{2}} \\
& \qquad \times 
 \sum_{l=0}^k [H_{l,k}^\mu]^{-1} 
  \left[  S_{k-2l,0}(\xi_\nu) P_{l,0}^k(x,y)+
     S_{k-2l,1}(\xi_\nu) P_{l,1}^k(x,y)\right] \\         
& = \frac{1}{n+1} \sum_{\nu=0}^{n} 
         a_\mu \int_{-1}^1  \CR_{\xi_\nu}(f;t) C_k^{\mu+1/2}(t) dt 
          \frac{k+\mu+\frac{1}{2}}{\mu+\frac{1}{2}} D_k^{\mu+1/2}(\xi_\mu;x,y),
\end{align*}
since $P_{l,\varepsilon}^k(\cos \xi_\nu,\sin\xi_\nu) = 
H_{l,k}^\mu S_{k-2l,\varepsilon}(\xi_\nu)$ and $\lambda_{l,k}^\mu = 
[H_{l,k}^\mu]^{-2}$. This proves the first identity.

We now prove the second equation \eqref{eq:proj2}. Using the fact that 
$\xi_{n+\nu+1}  = \xi_\nu + \pi$, 
$$
\cos (k-2l)\xi_{n+\nu+1} = (-1)^k \cos(k-2l)\xi_\nu, \quad 
 \sin (k-2l) \xi_{n+\nu+1} =  (-1)^k \sin (k-2l) \xi_\nu,
$$
we conclude that $D_k^{\mu+1/2}(\xi_{\nu};x,y)= (-1)^k C_k^{\mu+1/2}
(\xi_{n+1+\nu}; x,y).$
Hence, using the fact that $\CR_{\xi_\nu + \pi}^\mu(f;t) = 
\CR_{\xi_\nu }^\mu(f;-t)$, we conclude that 
$$
 \proj_{k}^\mu f(x,y)  = \frac{1}{n+1}  \sum_{\nu=0}^{n} 
  a_\mu \int_{-1}^1 \CR_{\xi_{n+1+\nu}}^\mu(f;t) 
         D_k^{\mu+1/2}(\xi_{n+1+\nu},t; x,y) dt.  
$$
Adding this equation and the first equation of \eqref{eq:proj} and dividing
the result by 2, we then have \eqref{eq:proj2}.  
\end{proof}

In the case of $\mu = 1/2$, it is easy to see that $\lambda_{l,k}^{\frac12} =
1/( k+1)$,  independent of $l$. Hence, for $\mu =1/2$, \eqref{eq:OPsum} 
shows that 
$$
  D_k^{\mu+1/2} (\xi_\nu; x,y) = \frac{1}{k+1} (k+1) C_k^{1} (\xi_\nu; x,y)
       =U_k (\xi_\nu; x,y),
 $$
and the formulas \eqref{eq:proj} and \eqref{eq:proj2} are of particular simple
form. This case was studied in \cite{X05}. 

The two expressions of $\proj_k f$ look similar but are different in an 
important point. The first expression consists of Radon projections in 
equally spaced directions along half of the the circumference of the circle, 
while the second expression uses Radon projections in equally spaced directions
along the entire circumference of the circle. This distinction is meaningful
for reconstruction algorithms for Radon data. 

If $n$ is even, then we can use Proposition \ref{prop:keyEven} instead of 
Proposition \ref{prop:key} in the proof. The result is another identity that 
uses Radon projections over equally spaced angles in $[0,2\pi]$. Let 
$$
     \phi_\nu = \frac{2 \nu \pi}{2m+1}, \qquad 0 \le \nu \le 2m.
$$

\begin{thm}  \label{thm:proj-sumEven}
For $m \ge 0$ and $k \le 2m$,  the operator $\proj_k^\mu f$ can be 
written as
\begin{align} \label{eq:projEven} 
\proj_{k}^\mu f(x,y) = \frac{1}{2m+1}  \sum_{\nu=0}^{2m} 
    a_\mu \int_{-1}^1 \CR_{\phi_\nu}^\mu(f;t) D_k^{\mu+1/2}(\phi_\nu,t; x,y) dt
\end{align}  
\end{thm}

This expression of $\proj_k f$ is not a special case of \eqref{eq:proj2}, even 
though both uses equally spaced angles. In fact, setting $n = 2m$ shows that 
\eqref{eq:proj2} uses exactly twice as Radon projections in equally spaced
directions. For $\mu =1/2$ the identity \eqref{eq:projEven} has appeared in 
\cite{X05}. The equation \eqref{eq:projEven} can be deduced from 
\eqref{eq:proj} as follows: using the fact that $\CR_{\phi+ \pi} f(t) = 
\CR_\phi(f;-t)$ and changing variable $t \mapsto -t$ in the integral 
whenever $\phi = \xi_{2\nu-1}$ in \eqref{eq:proj}, then making use of the 
equations in \eqref{eq:CosSin} and the fact that the Gegenbauer polynomial 
is symmetric.

Let $S_n^\mu f$ denote the $n$-th partial sum of the expansion 
\eqref{eq:fourier}; that is,
$$
S_n^\mu(f; x,y) =  \sum_{k=0}^n \proj_k^\mu f (x,y).
$$
The operator $S_n^\mu$ is a projection operator from $L^2(W_\mu;B^2)$ 
onto $\Pi_n^2$.  An immediate consequence of Theorem \ref{thm:proj-sum}
is the following corollary:

%%%%%%%%%%%%%%%%%%%%%%%%
%%%  Corollary 2.6  %%%%
%%%%%%%%%%%%%%%%%%%%%%%%

\begin{cor}  \label{thm:partial-sum}
For $n \ge 0$, the partial sum operator $S_{n}^\mu f$ can be written as
\begin{align} \label{eq:S2m} 
S_{n}^\mu(f;x,y) & = \frac{1}{n+1} \sum_{\nu=0}^{n} a_\mu 
  \int_{-1}^1 \CR_{\xi_\nu}^\mu (f;t) \Phi_n^\mu (\xi_\nu, t; x,y) dt \\
  & =  \frac{1}{2n+2}\sum_{\nu=0}^{2n+1} a_\mu 
  \int_{-1}^1 \CR_{\xi_\nu}^\mu (f;t) \Phi_n^\mu (\xi_\nu, t; x,y) dt 
\notag
\end{align}  
where 
\begin{equation} \label{eq:Phi} 
 \Phi_n^\mu (\xi,t; x,y) = 
    \sum_{k=0}^{n} \frac{k+\mu+1/2}{\mu+1/2} 
     C_k^{\mu+1/2}(t) D_k^{\mu+1/2}(\xi;x,y).
\end{equation} 
\end{cor}

Likewise, an immediate consequence of Theorem \ref{thm:proj-sumEven}
is the following corollary:

\begin{cor}  \label{thm:partial-sumEven}
For $m \ge 0$, the partial sum operator $S_{2m}^\mu f$ can be written as
\begin{equation} \label{eq:S2mEven} 
S_{2m}^\mu(f;x,y) = \frac{1}{2m+1}  \sum_{\nu=0}^{2m} a_\mu 
  \int_{-1}^1 \CR_{\phi_\nu}^\mu (f;t) \Phi_{2m}^\mu (\phi_\nu,t; x,y) dt. 
\end{equation}  
\end{cor}

\subsection{Discretization and reconstruction algorithm} 

The equation \eqref{eq:S2m} expresses the partial sum of the Fourier 
orthogonal expansion as the integrals of attenuated Radon projections in
the equally spaced directions. In order to derive an algorithm that uses 
only values of attenuated Radon projections on a set of finite line 
segments, we approximate the integrals by a quadrature formula. If $f$ 
is a polynomial then $\CR_\phi^\mu(f;t)/(1-t^2)^\mu$ is a polynomial of the 
same degree by Lemma \ref{lem:2.4}, which shows that we should use a 
quadrature 
formula with respect to the weight function $(1-t^2)^\mu$; that is, 
$$%\begin{equation} \label{eq:quadrature}
     \int_{-1}^1 g(t) (1-t^2)^\mu dt \approx  
        \sum_{j=1}^{N}  \lambda_{j}  g(t_{j}) 
$$%\end{equation}                                    
where $t_{j}$ are real numbers and $\lambda_{j}$ are chosen so that 
the quadrature produces exact values of the integrals for polynomials 
of degree at least $M$. Such a quadrature is said to be of $N$ points 
and of precision $M$.  A Gaussian quadrature of $N$ points has the 
highest precision $M = 2N-1$ among all quadrature formulas of $N$
points. 

For our purpose we are interested in quadrature formulas of precision $2n$
that uses $n+1$ points.  A class of such formulas is given in the following 
proposition, which is based on the zeros of the quasi-orthogonal
polynomial $C_{n+1}^{\mu+1/2}(t) + a C_n^{\mu+1/2}(t)$, where 
$a$ is a real number \cite{Sz}. For certain range of $a$, such a polynomial
has $n+1$ real distinct zeros in the interval $[-1,1]$. 

\begin{prop} \label{prop:gaussian}
Let $t_{j,n}$, $0 \le j \le n$, be the distinct zeros of a quasi-orthogonal
polynomial $C_{n+1}^{\mu+1/2}(t) + a C_n^{\mu+1/2}(t)$. Then there
are positive numbers $\lambda_{j,n}$ such that the quadrature 
\begin{equation} \label{eq:gaussian}
     \int_{-1}^1 g(t) (1-t^2)^\mu dt \approx  
        \sum_{j=0}^{n}  \lambda_{j,n}  g(t_{j,n}) : = \CI^\mu_{n}(g) 
\end{equation}                                   
has precision $2n$ if $a \ne 0$. If $a = 0$ then the quadrature has 
precision $2n+1$. 
\end{prop}

Using an appropriate quadrature on the integrals in \eqref{eq:S2m} we 
obtain a reconstruction algorithm for the attenuated Radon data. We state
such an algorithm only in the case of the quadrature formula in 
\eqref{eq:gaussian}. 
 
%%%%%%%%%%%%%%%%%
%%%%  Algorithm 2D-G %%%%
%%%%%%%%%%%%%%%%%

\begin{algo} \label{algo:2D-G} 
Let $\mu \ge 0$ and $n \ge 0$. Let $t_{j,n}$ and $\lambda_{j,n}$ 
be as in \eqref{eq:gaussian}. For $(x,y) \in B^2$ define  
\begin{equation} \label{eq:AlgoG} 
 \CA_{n}^\mu(f; x,y) = \sum_{\nu=0}^{n} \sum_{j=0}^{n}
    \CR^\mu_{\xi_\nu}(f; t_{j,n}) T_{j,\nu}^\mu(x,y),
\end{equation} 
where 
$$
T_{j,\nu}^\mu(x,y) = \frac{a_\mu \lambda_{j,n}}{n+1} 
            (1-t_{j,n}^2)^{-\mu}\Phi_n^\mu(\xi_\nu,t_{j,n}; x,y).
$$
\end{algo}
 
For a given $f$, the approximation process $\CA_{n}^\mu f$ uses 
attenuated Radon data 
$$
 \left \{\CR_{\xi_\nu}^\mu (f;  t_{j,n}): 0 \le \nu \le n, 
     \quad 0 \le j \le n\right\}
$$
of $f$. The data consist of Radon projections on $n+1$ equally spaced 
directions (specified by $\xi_\nu$) along the circumference of a half 
circle and there are $n+1$ parallel lines (specified by $t_{j,n}$) in each 
direction. The algorithm produces a polynomial $\CA_{n}^\mu f$ which is an 
approximation to $f$. In the case of $\mu =1/2$ the algorithm 
\eqref{algo:2D-G} appeared earlier in \cite{BO}; the connection to the 
orthogonal partial sums, however, was neither established nor used there. 

\begin{thm}
The operator $\CA_{n}^\mu$ is a projection operator on $\Pi_{n}^2$. In 
other words, $\CA_{n}^\mu f \in \Pi_{n}^2$ and $\CA_{n}^\mu P  = P$
for $P \in \Pi_{n}^2$.
\end{thm}

\begin{proof}
The function $\Phi^\mu(\xi_\nu,t_{j,n};x,y)$ is evidently an element in 
$\Pi_n^2$. It follows immediately that $\CA_n^\mu f \in \Pi_n^2$. By 
definition, $S_n^\mu $ is a projection operator on $\Pi_{n}^2$. The 
operator $\CA_n^\mu f$ is obtained from $S_n^\mu f$ by applying the 
quadrature \eqref{eq:gaussian}, exactly for polynomials in $\Pi^2_{2n}$, 
on $(1-t^2)^{-\mu} \CR_{\xi_\nu}^\mu (f; t)\Phi_n^\mu(\xi_\nu,t;\cdot)$, 
which is a polynomial of degree $2n$ in $t$ variable by Lemma \ref{lem:2.4}
and \eqref{eq:Phi} whenever $f \in \Pi_{n}^2$. Hence, the quadrature
\eqref{eq:gaussian} is exact. Thus, $\CA_{n}^\mu f = S_n^\mu f = f$ 
if $f \in \Pi_n^2$.
\end{proof}

Alternatively, we can use a quadrature formula of proper order on the 
second expression of \eqref{eq:S2m} to derive an algorithm that uses
Radon projections on $2n+2$ directions equally distributed along the
circumference of the entire circle. Instead of stating such an algorithm
we consider the case of $n=2m$ and use the expression \eqref{eq:S2mEven}.
This leads to an algorithm that sums over $2m+1$ angles that are equally 
spaced over $[0,2\pi]$, as we shall discuss in the following subsection.

\subsection{Reconstruction algorithm using attenuated Radon projections}

For practical applications in CT, the discretization described in 
Algorithm 2.13 needs to be further specified or simplified. In fact,  
one has to take into consideration what scan geometry is used in practice. 
For example, the zeros of quasi orthogonal polynomials will not be coincide
with the discrete measurement of the attenuated Radon projections in the
usual scan geometry. If these points were used, then it would be necessary
to introduce an interpolation process, which would introduce new errors. 
As an alternative, we suggest to use a different discretization, which
amounts to use a different quadrature formula.  

For the ordinary Radon projections ($\mu =1/2$), Gaussian quadrature formulas
for the weight function $\sqrt{1-x^2}$ are used for the integrals in 
\eqref{eq:S2mEven} to generate algorithms. For practical implementation in 
CT, the quadrature formula
\begin{equation}\label{eq:gaussianU}
\frac{1}{\pi} \int_{-1}^1 f(t) \frac{dt}{\sqrt{1-t^2}} =
\frac{1}{n+1} \sum_{j=0}^{n} f \left(\cos \tfrac{(2j+1)\pi}{2n+2}\right),
\end{equation}
based on zeros of $T_{n+1}(x) = \cos (n+1) \theta$, $x = \cos \theta$, is
used \cite{XTC}. The reason for such a choice lies in the scanning geometry 
of the input data. It turns out that, for $n = 2m$, such a choice allows 
us to adopt fan beam geometry and use it as parallel geometry in a 
straightforward way. 

It is possible to use the quadrature formula \eqref{eq:gaussianU} for 
attenuated Radon transforms $\CR_\phi^\mu (f;t)$, especially when $\mu$ is 
a half integer. The resulted $\CA_{2m}$ will no longer be a projection 
operator, but it still reproduces polynomials of degree slightly less than 
$n$ when $\mu$ is a half integer. 

%%%%%%%%%%%%%%%%%
%%%%%  Algorithm N/2 %%%%
%%%%%%%%%%%%%%%%%

\begin{algo} \label{algo:N/2} 
%Let $\mu$ be a half integer,  $\mu +1/2 \in \NN$. 
For $m \ge 0$, $(x,y) \in B^2$, 
\begin{equation} \label{eq:AlgoN/2} 
 \CA_{2m}^\mu(f; x,y) = \sum_{\nu=0}^{2m} \sum_{j=0}^{2m}
 \CR_{\phi_\nu}^\mu \left(f;\cos \psi_j \right) T_{j,\nu}^\mu(x,y),
\end{equation} 
where 
$$
T_{j,\nu}^\mu(x,y) =  \frac{\mu+1/2}{(2m+1)^2} 
     \sin \psi_j \Phi_{2m}^\mu(\phi_\nu, \cos \psi_j; x,y), \qquad
      \psi_j = \frac{(2j+1) \pi}{4m+2}.
$$
\end{algo}

The constant $\mu + 1/2$ in $T_{j,\nu}^\mu$ comes from the fact that 
$a_\mu = (\mu+1/2)/\pi$. 

This algorithm provides an approximation for the reconstruction of a 
function $f(x,y)$ from a set of attenuated Radon projections 
$$
 \left \{ \CR_{\phi_\nu}^\mu (f; \cos \psi_j),  
    \quad 0 \le \nu \le 2m, \quad 1\le j \le 2m  \right \}.
$$
The set $\{\phi_\nu: 0 \le \nu \le 2m \}$ consists of equally spaced angles 
along the circumference  of the disk.  For $\mu =1/2$ it has appeared in 
\cite{X05}. The advantage of this algorithm lies in the fact that it can be 
used with attenuated Radon data obtained from the fan beam geometry 
directly, see the discussion in \cite{XTC}.  The operator, however, 
reproduces polynomials up to a lower degree. 

\begin{thm}
Let $\mu$ be a half integer,  $\mu +1/2 \in \NN$. Then the operator 
$\CA_{2m}^\mu$ in Algorithm \ref{algo:N/2} preserves polynomials of degree 
$2m-2\mu$; that is,  $\CA_{2m}^\mu P  = P$ for $P \in \Pi_{2m- 2\mu}^2$.
\end{thm}

\begin{proof}
The algorithm is obtained by using  the Gaussian quadrature formula 
\eqref{eq:gaussianU} to discretize the integrals in \eqref{eq:S2mEven}, 
that is,
\begin{align*}
  \int_{-1}^1 \CR_{\phi_\nu}^\mu (f; t) C_k^{\mu+1/2}(t) dt  & =
     \int_{-1}^1 \frac{\CR_{\phi_\nu}^\mu(f; t)}{\sqrt{1-t^2}}
        C_k^{\mu+1/2}(t) \sqrt{1-t^2} dt \\
      & \approx \frac{\pi}{2m+1} \sum_{k=0}^{2m} \sin\psi_j 
        \CR_{\phi_\nu}^\mu(f;\cos \psi_j) C_k^{\mu+1/2}(\cos \psi_j). 
 \end{align*}
If $f\in \Pi_{2m-2\mu}^2$ then using the fact that 
$\CR_{\phi}^\mu(f;t)/(1-t^2)^\mu$ is a polynomial of degree $2m-2\mu$, 
the assumption that $\mu$ is a half integer shows that 
$$
\CR_{\phi_\nu}(f; t)/\sqrt{1-t^2} =  (1-t^2)^{\mu-1/2} 
\CR_{\phi_\nu}(f; t)/(1-t^2)^{\mu}
$$ 
is a polynomial of $2\mu-1+2m - 2\mu = 2m -1$. Since 
$\Phi_{2m}^\mu(\xi_\nu,t; \cdot)$ is of degree $2m$ and the quadrature
\eqref{eq:gaussianU} is of precision $4m-1$, the discretization becomes 
exact in this case and we conclude that $\CA_{2m}^\mu f =f$ if $f \in 
\Pi_{2m-2\mu}^2$.
\end{proof}

Let $C(B^2)$ denote the space of continuous function on $B^2$ with 
the uniform norm $\|\cdot\|_\infty$ and let  $\|\CA_{n}^\mu\|$ denote 
the operator norm of $\CA_{n}^\mu$ under the uniform norm.  By 
$A \sim B$ we mean that there are two constants $c_1$ and $c_2$ 
such that $c_1 A \le B \le c_2 A$.  Evidently the convergence of the 
algorithm depends on $\|\CA_{n}^\mu\|$. In fact, since $\CA_{n}^\mu$ 
in Algorithm \ref{algo:2D-G} preserves $\Pi_n$, it is easy to see that
$$
 \|f - \CA_{n}^\mu f\| \le c_f \left( 1 + \|\CA_{n}^\mu\| \right) E_n(f)
$$ 
where $E_n(f): = \inf \{\|f-P\|: P \in \Pi_n^2\}$ is the error of the best 
approximation of $f$ by polynomials on $B^2$.  If $f$ has $r$-th order
continuous derivatives, then $E_n(f) \le c_f n^{-r}$, in which $c_f$ 
depends on the norm of the $r$-th derivatives of $f$.  The same applies 
to $\CA_{2m}^\mu$ in Algorithm \ref{algo:N/2}, which preserves 
$\Pi_{2m -2\mu}$. Using the formula in \eqref{eq:Radon}, the proof of 
Proposition 5.1 of \cite{X05} gives the following formula of the norm of 
$\CA_{2m}^\mu$ in Algorithm \ref{algo:N/2}: 

\begin{prop}
The operator norm $\|\CA_{2m}^\mu\|$ of $C(B^2)$ to $C(B^2)$ is given by
$$
  \|\CA_{2m}^\mu \| = \max_{(x,y) \in B^2} \Lambda_m (x,y), \qquad
  \Lambda_{m}(x,y) := \sum_{\nu=0}^{2m} \sum_{j=0}^{2m}
    (\sin \theta_{j,m})^\mu |T_{j,\nu}^\mu (x,y)|.   
$$
\end{prop}

As $m \to \infty$, the norm growth in an essentially polynomial order 
of $m$. Hence, the algorithm converges uniformly if $f$ is sufficiently 
smooth. To estimate the exact order of $\CA_{2m}^\mu$ is difficult. 
In the case of $\mu =1/2$, it is carried out in \cite{X05} and the 
order is $\|\CA_{2m}\| \sim  m\log (m+1)$. Based on this fact, we conjecture
that the operator norm of $\CA_{2m}^\mu$ is of the the order 
$$
     \|\CA_{2m}^\mu \| \sim m^{\mu+1/2} \log (m+1), 
            \qquad\hbox{as $ m \to \infty$},
$$
which is only slightly worse than the norm $\|S_{n}^\mu\| \sim 
n^{\mu+1/2}$  (\cite{X01}). If the conjecture holds, then the algorithm
will converges uniformly for smooth $f \in C^r(B^2)$ with $r > \mu + 1/2$.
In most applications, however, the function or image could have jumps;
that is, there is not even continuity. The numerical tests in the case
of ordinary Radon data shows that the algorithm is stable and yields 
fairly accurate results even when the data is highly singular (\cite{XTC}).
See also the example given in the following subsection.

\subsection{Numerical Example}

For the numerical examples we use Algorithm \ref{algo:N/2}, for which
the scan geometry is easy to implement. The data required are 
$g_{j,\nu}:=\CR^\mu_{\phi_\nu}(f; \cos \psi_j)$, where $\phi_\nu = 
2\nu \pi/(2m+1)$ stands for the $2m+1$ views equally spaced along the 
circumference of the region to be reconstructed and $\psi_j = (2j+1)/(4m+2)$
means that the x-rays in each view is distributed according to the zeros
of the Chebyshev polynomial $T_{2m+1}$. In this case the fan data can be 
resorted into parallel data (\cite{XTC}). 

We reconstruct a simple analytical phantom defined by the function 
$$
  f(x,y) = \begin{cases} 1 & \hbox{if  $0.9 \le r \le 1$ or 
         $0 \le r \le 0.1$} \\
             0 &  \hbox{if $0.1 < r < 0.9$},
           \end{cases} 
$$
where $r = \sqrt{x^2 + y^2}$, on the unit disk. This phantom contains 
strong singularity along the circles $r = 0.9$ and $r = 0.1$. The rotationally
invariant nature of the function allows certain simplification of the 
algorithm. 

For the reconstruction, we choose three values of the parameter $\mu$,
$\mu = 0, 1/2, 3/2$. The case $\mu =0$ means the ordinary Radon transform.
The case $\mu =0$ means that the Radon transform is attenuated by the
weight function $(1-x^2-y^2)^{-1/2}$, which is infinity at the boundary
of the disk. The case $\mu =3/2$ means that the Radon transform is 
attenuated by the weight function $1-x^2-y^2$, which is zero at the 
boundary. In each case, the Radon data are computed analytically. 

For each of the three values of $\mu$, we use Algorithm \ref{algo:N/2} for 
the reconstruction with a moderate $m =100$. The reconstructed image is 
evaluated on a $300 \times 300$ grid. The result is shown in Figure 1 below. 

\medskip

%\centerline{
%\ifpdf
%\includegraphics[width = 4cm]{rec_mu_0_m100.pdf} \hskip .05in
%\includegraphics[width = 4cm]{rec_mu_1_2_m100.pdf} \hskip .05in
%\includegraphics[width = 4cm]{rec_mu_3_2_m100.pdf}
%\else
%\includegraphics[width = 4cm]{rec_mu_0_m100.epsi} \hskip .05in
%\includegraphics[width = 4cm]{rec_mu_1_2_m100.epsi} \hskip .05in
%\includegraphics[width = 4cm]{rec_mu_3_2_m100.epsi}
%\fi}
\centerline{\epsfxsize=4cm \epsffile{rec_mu_0_m100.epsi} \hskip .05in
\epsfxsize=4cm \epsffile{rec_mu_1_2_m100.epsi} \hskip .05in
\epsfxsize=4cm \epsffile{rec_mu_3_2_m100.epsi} }
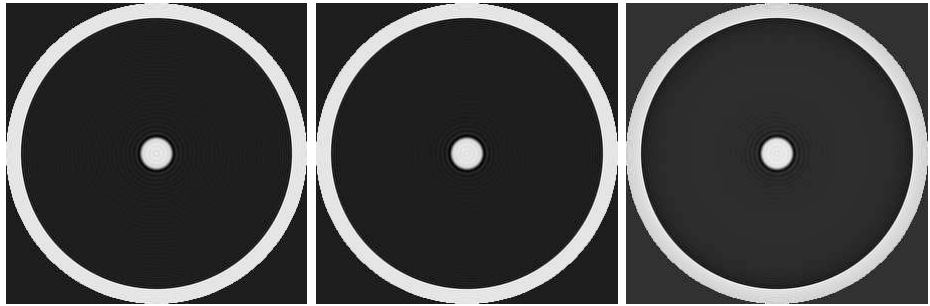
\centerline{ Figure 1. From left to right, $\mu = 0,1/2,3/2$.}

\medskip

These images show that the function is reconstructed rather faithfully
in each of the three cases, even though the function has strong singularity.
The case $\mu = 1/2$ has been tested extensively and compared with FBP 
algorithm (\cite{XTC,XTC2}). The above is our first attempt to test the 
algorithm for attenuated Radon transforms. 

%%%%%%%%%%%%%%%%%%%%%%% 
%%%%%% Section 3  %%%%%
%%%%%%%%%%%%%%%%%%%%%%%

\section{Reconstruction and Approximation on the unit sphere}
\setcounter{equation}{0}

It is known that orthogonal polynomials on the unit ball and on the unit
sphere are closely related (\cite{X98}). Since the approximation and the
reconstruction in the previous section are based on orthogonal expansions
on the unit disk, the relation suggests analogous results on the unit 
sphere $S^2 =\{(x,y,z): x^2+y^2+z^2 =1\}$, which we explore in this section. 

On the sphere we consider the attenuated spherical transform defined by
$$
Q^\mu f(\zeta;t) = \int_{\langle \xb, \zeta \rangle = t} 
f(\xb) |x_3|^{2\mu} 
    d \omega,
$$
where $\xb = (x_1,x_2,x_3) \in S^2$, $\zeta \in \RR^3$ and $\xi \ne 0$, 
and $d\omega$ is the measure on the subset $\{x \in S^2: 
\langle \xb, \zeta \rangle = t\}$ which is the circle on the sphere. 
When $\mu = 0$, this is the usual spherical transform \eqref{eq:1.4}, see 
for example, \cite[p. 33]{NW}. We will mainly work with the case that 
$\zeta_3 =0$. We say that a function is even in $x_3$ if $f(x_1,x_2,x_3) =
 f(x_1,x_2,-x_3)$. 

\begin{prop}
Let $f$ be even in $x_3$. If $\zeta = (\cos\theta, \sin\theta, 0)$, then
\begin{equation} \label{eq:Q-mu}
    Q^\mu f(\zeta;t) = \CR_\theta^\mu(F;t), \qquad 
           F(x_1,x_2) = f\left(x_1,x_2,\sqrt{1-x_1^2-x_2^2} \right)
\end{equation}
\end{prop}

\begin{proof}
Since $f$ is even in $x_3$ we have  $f(\xb) = F(x_1,x_2)$ for $\xb \in S^2$. 
The definition of $\zeta$ shows that $\langle \xb, \zeta\rangle = 
x_1 \cos \theta + x_2 \sin \theta = I(\theta,t)$. In terms of $x_1$ and
$x_2$, $d\omega =  dx_1 dx_2/\sqrt{1-x_1^2-x_2^2}$. Thus, 
$$
Q^\mu f(\zeta;t) = \int_{x_1 \cos \theta + x_2 \sin\theta = t} F(x_1,x_2) 
  \left(1-x_1^2-x_2^2\right)^{\mu} \frac{dx_1  dx_2}{\sqrt{1-x_1^2-x_2^2}},
$$
which is precisely $\CR_\theta^\mu (F;t)$.
\end{proof}

Let $H_\mu(\xb) = |x_3|^{\mu}$. The space $L^2(H_\mu;S^2)$ has an
orthogonal decomposition 
\begin{equation} \label{eq:fourierS}
    L^2(H_\mu; S^2) = \sum_{k=0}^\infty \bigoplus \CH_k^\mu
\end{equation}
where the subspaces $\CH_k^\mu$ contains homogeneous polynomials 
of degree $k$ that are orthogonal to lower degree polynomials with respect 
to $H_\mu d\omega$ on $S^2$. For $\mu = 0$, $\CH_k^0$ is the space of 
ordinary spherical harmonics. Let 
$$
  \proj_{\CH_k^\mu}f: L^2(H_\mu;S^2) \mapsto \CH_k^\mu 
$$
be the orthogonal projection from $L^2(H_\mu;S^2)$ onto $\CH_k^\mu$. The space
$\CH_k^\mu$ is closely related to the space $\CV_k^2(W_\mu)$ discussed 
in the previous section (\cite{X98}). For our purpose, we only need the 
following relation on the orthogonal projections: if $f$ is even in $x_3$ then 
\begin{equation}\label{eq:projections}
\proj_{\CH_n^\mu}  f(\xb) = \proj_n^\mu F(x_1,x_2),
\end{equation}
where $F$ is the function defined in \eqref{eq:Q-mu}. This relation, together
with \eqref{eq:Q-mu}, allows us to express the projection operator on the
sphere in terms of spherical transforms. Using these relations and Theorem
\ref{thm:proj-sum} we obtain the following result: 
    
%%%%%%%%%%%%%%%%%
%%%%%  Theorem 3.2  %%%%
%%%%%%%%%%%%%%%%%

\begin{thm}  \label{thm:proj-sumS}
Let $f$ be even in $x_3$. For $n \ge 0$ and $k \le n$, the operator 
$\proj_{\CH_k^\mu}$ can be written as
\begin{align} \label{eq:projS} 
\proj_{\CH_k^\mu} f(\xb) = \frac{1}{n+1}  \sum_{\nu=0}^{n} 
  a_\mu \int_{-1}^1 Q^\mu f(\zeta_\nu; t) D_k^{\mu+1/2}(\xi_\nu,t; x_1,x_2) dt 
 \end{align}  
where $\xi_\nu = \tfrac{\nu\pi}{n+1}$, 
$\zeta_\nu = (\cos \xi_\nu, \sin \xi_\nu,0)$, and $D_k^{\mu+1/2}(\xi,t;x,y)$
is defined in \eqref{eq:C-t-xi}.  
\end{thm}

Let $Y_n^\mu f$ denote the $n$-th partial sum of the expansion 
\eqref{eq:fourierS}; that is, 
$$
  Y_n^\mu(f; \xb) =  \sum_{k=0}^n \proj_{\CH_k^\mu} f (x_1,x_2).
$$
The operator $Y_n^\mu$ is a projection operator from $L^2(H_\mu;S^2)$ 
onto $\Pi_n(S^2)$, the restriction of $\Pi_n^3$ on $S^2$.  An immediate 
consequence of Theorem \ref{thm:proj-sumS} is the following:

%%%%%%%%%%%%%%%%%
%%%%%  Corollary 3.3  %%%%
%%%%%%%%%%%%%%%%%

\begin{cor}  \label{thm:partial-sumS}
Let $f$ be even in $x_3$. For $n \ge 0$, the partial sum operator 
$Y_{n}^\mu f$ can be written as 
\begin{align} \label{eq:S2mS} 
Y_{n}^\mu(f; \xb) & = \frac{1}{n+1}\sum_{\nu=0}^{n} a_\mu 
  \int_{-1}^1 Q^\mu f(\zeta_\nu; t) \Phi_n^\mu (\xi_\nu,t; x_1,x_2) dt \\
 & = \frac{1}{2n+2}\sum_{\nu=0}^{2n+1} a_\mu 
  \int_{-1}^1 Q^\mu f(\zeta_\nu; t) \Phi_n^\mu (\xi_\nu,t; x_1,x_2) dt \notag
\end{align}  
where $\Phi_n^\mu$ is the function defined in \eqref{eq:Phi}. 
\end{cor}

For $n = 2m$ we can also use Theorem \ref{thm:proj-sumEven} to 
derive an expression for $Y_{2m}^\mu (f)$, which leads to the following
corollary:

\begin{cor}  \label{thm:partial-sumSEven}
Let $f$ be even in $x_3$. For $m \ge 0$, the partial sum operator 
$Y_{2m}^\mu f$ can be written as
\begin{equation} \label{eq:S2mSEven} 
Y_{2m}^\mu(f; \xb) =  \frac{1}{2m+1} \sum_{\nu=0}^{2m} a_\mu 
  \int_{-1}^1 Q^\mu f(\zeta_\nu; t) \Phi_{2m}^\mu (\phi_\nu,t; x_1,x_2) dt
\end{equation}  
where $\phi_\nu = \tfrac{2\nu \pi}{2m+1}$, $\zeta_\nu = 
(\cos \phi_\nu, \sin\phi_\nu, 0)$, and  $\Phi_{2m}^\mu$ is the function 
defined in \eqref{eq:Phi}. 
\end{cor}

In the case of $\mu = 0$, the equations \eqref{eq:S2mS} and 
\eqref{eq:S2mSEven} are representations of the partial sums of ordinary 
spherical harmonic expansions, which are expressed in terms of the 
Legendre polynomial $P_k(t) = C_k^{1/2}(t)$. 

Just like the case of orthogonal expansions on the unit disk, we can use a 
quadrature formula to obtain a reconstruction  algorithm using 
spherical transforms. For example, using the quadrature formula 
with respect to $(1-t^2)^\mu$ in Proposition \ref{prop:gaussian} as in the 
case of Algorithm \ref{algo:2D-G}, we get the following result: 

%%%%%%%%%%%%%%%%%
%%%%  Algorithm 3DS-G %%%
%%%%%%%%%%%%%%%%%

\begin{algo} \label{algo:3DS-G} 
Let $f$ be even in $x_3$. Let $\mu \ge 0$. For $n \ge 0$, $\xb \in S^2$,
\begin{equation} \label{eq:AlgoS-G} 
 \CS_{n}^\mu(f; \xb) = \sum_{\nu=0}^{n} \sum_{j=0}^{n}
      Q^\mu f(\zeta_\nu; t_{j,n}) T_{j,\nu}^\mu(x_1,x_2),
\end{equation} 
where $t_{j,n}$  are as in the quadrature \eqref{eq:gaussian} and 
$T_{j,\nu}^\mu$ are  defined in Algorithm \ref{algo:2D-G}.
\end{algo}

This algorithm reconstructs a function $f(\xb)$ from a set of spherical 
transforms 
$$
 \left \{ Q^\mu f(\zeta_\nu;t_j), \quad \zeta_\nu =
    (\cos \xi_\nu,\sin \xi_\mu,0), 
    \quad 0 \le \nu \le 2m, \quad 1\le j \le 2m  \right \},
$$
which consists of integrals over a number of circles on the sphere. These
circles lie on planes that are parallel to the $x_3$-axis. The circles 
intersect the circumference of a disk 
perpendicular to the $x_3$-axis at equally spaced angles. The distance
 between parallel circles depends 
on the values of $t_{j,n}$.  In the case $\mu =0$, the algorithm provides
an approximation to the function based on ordinary spherical transforms.
The assumption that $f$ is even in $x_3$ implies that we can use the
algorithm to reconstruct a function defined on the upper hemisphere from 
spherical transforms that are integrals over half circles parallel to 
$x_3$ axis on the upper hemisphere. 

If $\mu$ is a half integer, we can also state an algorithm using the 
quadrature \eqref{eq:gaussianU}, as in Algorithm \ref{algo:N/2}, so that
$t_{j,n} = \cos  j \pi/ (2m+1)$. However, in the most interesting case of 
$\mu =0$,  we do not have such a somewhat simplified algorithm.  

%%%%%%%%%%%%%%%%%%%%%%%
%%%%%% Section 4  %%%%%
%%%%%%%%%%%%%%%%%%%%%%%

\section{Reconstruction and Approximation on the unit ball}
\setcounter{equation}{0}

In this section we consider reconstruction of functions on a unit 
ball $B^3$ in $\RR^3$ based on the attenuated Radon projections. 

\subsection{Radon projections and orthogonal polynomials}

We will work with attenuated Radon projections that are integrals on line
segments inside $B^3$ with respect to the weight function  
$$
   W_\mu(\xb) = (1-\|\xb\|^2)^{\mu-1/2}, \quad \xb = (x_1,x_2,x_3) \in B^3, 
      \quad \mu \ge 0.
$$
For our purpose, however, we will only consider those lines lying on the 
planes that are perpendicular to the $x_3$ axis. Let $x_3 =  w$ be such a 
plane. Its intersection with the unit ball $B^3$ is a disk 
$\{\xb: x_1^2+x_2^2 \le \sqrt{1-w^2}, x_3 = w\}$.  A line on this disk 
is given by the equation
$$
\ell:\quad x\cos\theta+ y\sin \theta = t \sqrt{1-w^2}, \qquad  -1 \le t \le 1.
$$
Let $I(\theta,w; t)$ denote the intersection of $\ell$ with $B^3$. The
attenuated Radon projection on such a line is then defined by
\begin{equation}\label{eq:Radon3D}
  \CR_{\theta}^\mu(f;t, w) :=   \int_{I(\theta,w;t)}  f(\xb) W_\mu(\xb)d\xb. 
\end{equation}
The case $\mu = 1/2$ again corresponds to the usual Radon projection. 

%%%%%%%%%%%%%%%%%%%%%%%%%%
%%%%    Lemma 4.1     %%%%
%%%%%%%%%%%%%%%%%%%%%%%%%%

\begin{lem} \label{lem:4.1}
For  $f \in L^1(W_\mu; B^3)$ and for a fixed $w \in [-1,1]$, define a 
function $g_w$ on $B^2$ by 
$$
  g_w(x,y) =  f \left( \sqrt{1-w^2}\,  x, \sqrt{1-w^2} \, y, w \right). 
$$
The X-ray transform \eqref{eq:Radon3D} is related to the 2D Radon 
transform \eqref{eq:R-mu} by 
\begin{equation}\label{eq:Radon3D-2D}
  \CR_{\theta}^\mu(f;t,w) = (1-w^2)^\mu \CR_{\theta}^\mu (g_w;t).
\end{equation}
\end{lem}

\begin{proof}
Since $I(\theta,w;t)$ can be represented by  
$$
 x_1 = \sqrt{1-w^2} ( t \cos \theta - s \sin \theta), \quad 
 x_2 = \sqrt{1-w^2} ( t \sin \theta + t \cos \theta), \quad x_3 =w 
$$
for $s \in [-\sqrt{1-t^2}, \sqrt{1-t^2}]$, which is a rotation around $x_3$ 
axis on the plane defined by $x_3 =w$, we have 
$$
      \CR_{\theta}(f;t,w) 
             =  (1-w^2)^\mu  \int_{-\sqrt{1-t^2}}^{\sqrt{1-t^2}}
                      g_{w} (t\cos \theta - s \sin \theta,  
                          t \sin \theta + s \cos \theta) W_\mu(s,t) ds.   
$$
The integral is precisely $\CR^\mu_\theta(g_{w};t)$ by \eqref{eq:Radon}. 
\end{proof}

Let $\CV_n^3(W_\mu)$ denote the space of orthogonal polynomials with 
respect to $W_\mu$ on $B^3$, which contains polynomials of degree
$n$ that are orthogonal to polynomials of lower degrees with respect to the
inner product 
$$
 \langle P, Q\rangle = a_{\mu,3} \int_{B^3} P(\xb) Q(\xb) W_\mu(x) d\xb, 
   \qquad
     a_{\mu,3} = \frac{\Gamma(\mu+2)}{ \pi^{3/2}\Gamma(\mu+1/2)},
$$
where $a_{\mu,3}$ is the normalization constant of $W_\mu$. We derive 
a basis for $\CV_n^3(W_\mu)$, making use of an orthogonal basis for 
$\CV_n^2(W_\mu)$. We note that the $W_\mu$ in these two notations 
are different, the first one is on $B^3$ and the second one is on $B^2$.  
We denote by  $\wt C_j^\lambda$ the orthonormal Gegenbauer
polynomial, which is equal to $C_n^\lambda / \sqrt{h_n}$ by
 \eqref{eq:Gegenbauer}. 

%%%%%%%%%%%%%%%%%
%%%%%  Proposition 4.2 %%%
%%%%%%%%%%%%%%%%%

\begin{prop}  \label{prop:Qlkj}
Let $\{P_j^k: 0 \le j \le k\}$ be an orthonormal basis for 
$\CV_k^2(W_\mu)$. Then the polynomials 
\begin{equation} \label{eq:Qlkj}
    Q_{l,k,j}(x,y,z) = h_k (1-z^2)^{k/2} P_j^k\left( \frac{x}{\sqrt{1-z^2}},
          \frac{y}{\sqrt{1-z^2}}\right) \wt C_{l-k}^{k + \mu+1}(z) 
\end{equation}
for $0 \le j \le k \le l$, where $h_k^2 = (\mu+2)_{k} / (\mu+3/2)_{k}$, form 
an orthonormal basis for $\CV_l^3(W_\mu)$. 
\end{prop}

\begin{proof}
From Lemma \ref{lem:2.3}, it is easy to see that $P_j^k$ is a sum of 
even powers of homogeneous polynomials when $k$ is even, and 
a sum of odd powers of homogeneous polynomials when $k$ is odd.
Thus, it follows that $Q_{l,k,j} \in \Pi_l^3$. 
Using the fact that $P_j^k$ is orthonormal, it follows from the integral 
relation 
\begin{equation} \label{eq:B3Integral}
\int_{B^3} f(\xb) d\xb =  \int_{-1}^1 \int_{B^2} 
  f\left(x_1 \sqrt{1-x_3^2}, x_2 \sqrt{1-x_3^2}, x_3\right)
      dx_1 dx_2 (1-x_3^2) dx_3 
\end{equation}
and the fact that $a_{\mu,3} = a_\mu c_{\mu+1/2}$, where $a_\mu$ is the
normalization of $W_\mu$ on $B^2$ and $c_\mu$ is defined in 
\eqref{eq:Gegenbauer}, that
\begin{align*}
 & a_{\mu,3} \int_{B^3} Q_{l,k,j}(\xb) Q_{l',k',j'}(\xb) W_\mu (\xb) d\xb \\
 & \quad =  h_k^2 c_{\mu+1/2} \int_{-1}^1 \wt C_{l-k}^{k+\mu+1}(t)
     \wt C_{l'-k}^{k+\mu+1}(t) (1-t^2)^{k+\mu+1/2} dt
  \delta_{k,k'}\delta_{j,j'} \\
       & \quad =  h_k^2 \frac{c_{\mu+1/2}}{c_{k+\mu+1/2}} 
             \delta_{l,l'} \delta_{k,k'}\delta_{j,j'}.
\end{align*}
It follows from the definition of $c_\mu$ that $c_{\mu+1/2}/c_{k+\mu+1/2}
 = (\mu+3/2)_k/(\mu+2)_k$, which completes the proof.
\end{proof} 

The attenuated Radon transforms of this basis can be computed explicitly.

\begin{prop} \label{prop:4.3}
Let $\mu \ge 0$ and let $Q_{l,k,j}$ be defined by \eqref{eq:Qlkj}. 
Then 
\begin{align} \label{eq:RQ3}
 & \frac{ \CR_{\phi}^\mu(Q_{l,k,j};t,w)}{(1-t^2)^\mu(1-w^2)^\mu} \\
& \qquad  = b_\mu \frac{C_k^{\mu+1/2}(t)}{C_k^{\mu+1/2}(1)}
         Q_{l,k,j}\left(\sqrt{1-w^2} \cos \phi, \sqrt{1-w^2} 
            \sin \phi, w\right).  \notag 
\end{align}
\end{prop}

\begin{proof}
By Lemma \ref{lem:4.1} and the definition of $Q_{l,k,j}$ we have 
$$
 \CR_{\phi}^\mu(Q_{l,k,j};t,w) = (1-w^2)^\mu \CR_\phi^\mu(g_w;t),
$$
where $g_w(x,y) = h_k P_j^k(x,y) (1-w^2)^{k/2}\wt C_{l-k}^{k+\mu+1}(w)$. 
By Lemma \ref{lem:2.4}, it follows that 
\begin{align*}
  \CR_\phi^\mu(g_w;t) & = h_k (1-w^2)^{k/2} \wt C_{l-k}^{k+\mu+1}(w) 
           \CR_\phi^\mu(P_j^k ;t) \\
        & = b_\mu h_k (1-w^2)^{k/2} \wt C_{l-k}^{k+\mu+1}(w) (1-t^2)^\mu
    \frac{C_k^{\mu+1/2}(t)}{C_k^{\mu+1/2}(1)} P_j^k(\cos \phi,\sin\phi) \\
             & = b_\mu (1-t^2)^\mu
    \frac{C_k^{\mu+1/2}(t)}{C_k^{\mu+1/2}(1)} 
      Q_{l,k,j}\left(\sqrt{1-w^2}\cos \phi,\sqrt{1-w^2}\sin\phi,w\right)
\end{align*}
by the definition of $Q_{l,k,j}$. Putting these equations together completes
the proof.
\end{proof}

Let $\proj_{l,3}^\mu $ denote the projection operator from $L^2(W_\mu; B^3)$ 
onto the space $\CV_l^3(W_\mu)$. Again we have the decomposition
\begin{equation} \label{eq:expansion}
   L^2(W_\mu; B^3) = \sum_{k=0}^\infty \bigoplus \CV_k^3(W_\mu): 
     \qquad f = \sum_{k=0}^\infty \proj_{k,3}^\mu f.
\end{equation}

%%%%%%%%%%%%%%%%%
%%%%%  Proposition 4.4 %%%
%%%%%%%%%%%%%%%%%

\begin{prop} \label{prop:4.4} 
For $n \ge 0$ and $0 \le l \le n$, 
\begin{align} \label{eq:ProjB3}
\proj_{l,3}^\mu f(\xb) & = 
 \frac{1}{n+1}\sum_{\nu =0}^{n} \int_{-1}^1\int_{-1}^1
 \CR_{\xi_\nu}^\mu (f;t,w) G_{l} (\xi_\nu, t, w; \xb) dt \sqrt{1-w^2} dw \\
& = \frac{1}{2n+2}\sum_{\nu =0}^{2n+1} \int_{-1}^1\int_{-1}^1
 \CR_{\xi_\nu}^\mu (f;t,w) G_{l} (\xi_\nu, t, w; \xb) dt \sqrt{1-w^2} dw
 \notag
\end{align}
where 
\begin{align*}
   G_{l}(\xi,t,w;\xb) = & a_{\mu,3} \sum_{k=0}^l h_k^2
       D_k^{\mu+1/2}\left(\xi, t; \tfrac{x_1}{\sqrt{1-x_3^2}}, 
            \tfrac{x_2}{\sqrt{1-x_3^2}}\right) \\  
 & \times (1-w^2)^{k/2} (1-x_3^2)^{k/2}  \wt C_{l-k}^{k+\mu+1}(w) 
          \wt C_{l-k}^{k+\mu+1}(x_3). 
\end{align*}
\end{prop}

\begin{proof}
The projection operator has an integral expression just as that of 
\eqref{eq:Pintegral}. Furthermore, the kernel function $P(W_\mu;\xb,\yb)$
can be written as a sum of an orthonormal basis. In particular, 
$$
     \proj_{l,3}^\mu f(\xb) = \sum_{k=0}^l 
          \sum_{j=0}^k \wh f_{l,k,j} Q_{l,k,j}(\xb),
$$
where $Q_{l,k,j}$ is the orthonormal basis for $\CV_l^3(W_\mu)$ defined in
\eqref{eq:Qlkj} and 
$$
 \wh f_{l,k,j} = a_{\mu,3} \int_{B^3} f(\yb) Q_{l,k,j}(\yb) W_\mu(\yb) d\yb.  
$$
Using \eqref{eq:B3Integral}, the definition of $Q_{k,l,j}$, and the
fact that $a_{\mu,3} = a_\mu c_{\mu+1/2}$, we have 
\begin{align*}
   \wh f_{l,k,j} = c_{\mu+1/2} \int_{-1}^1 & \left[ a_{\mu} \int_{B^2} 
       g_w(u,v) P_j^k(u,v) W_\mu(u,v) du dv \right] \\
        & \times  h_k \wt C_{l-k}^{k+\mu+1}(w)(1-w^2)^{k/2+\mu+1/2} dw,  
\end{align*}
where $g_w$ is defined as in Lemma \ref{lem:4.1}. Hence, it follows from
\eqref{eq:Pintegral} and \eqref{eq:Pkernel} that 
\begin{align*}
 \proj_{l,3}f(\xb) &= \sum_{k=0}^l h_k^2 \wt C_{l-k}^{k+\mu+1}(x_3)
     (1-x_3^2)^{k/2} c_{\mu+1/2}\\
 & \times \int_{-1}^1 \proj_k^\mu g_w\left(\tfrac{x_1}{\sqrt{1-x_3^2}},
  \tfrac{x_2}{\sqrt{1-x_3^2}}  \right) \wt C_{l-k}^{k+\mu+1}(w) 
        (1-w^2)^{(k+1)/2+\mu} d w.
\end{align*}
The identity \eqref{eq:ProjB3} follows from the above equation upon using
\eqref{eq:proj} and \eqref{eq:Radon3D-2D}.
\end{proof}

Let us denote by $S_{n,3}^\mu f$ the $n$-th partial sum of the orthogonal
expansion \eqref{eq:expansion},
$$
  S_{n,3}^\mu f(\xb) = \sum_{l=0}^n \proj_{l,3}^\mu f(\xb).
$$
As an immediate consequence of Proposition \ref{prop:4.4} we have 

\begin{cor} 
For $n \ge 0$, 
\begin{align}\label{eq:S-2m3}
S_{n,3}^\mu f(\xb) = \frac{1}{n+1}\sum_{\nu = 0}^{n} \int_{-1}^1\int_{-1}^1
 \CR_{\xi_\nu}^\mu(f;t,w)\Phi_{n}^\mu (\xi_\nu,t,w;\xb)dt\sqrt{1-w^2}dw \\
= \frac{1}{2n+2}\sum_{\nu = 0}^{2n+1} \int_{-1}^1\int_{-1}^1
 \CR_{\xi_\nu}^\mu(f;t,w)\Phi_{n}^\mu (\xi_\nu,t,w;\xb)dt\sqrt{1-w^2}dw 
 \notag
\end{align}
where 
$$
 \Phi_{n}^\mu(\xi,t,w;\xb) = \sum_{l=0}^{n} G_{l}(\xi,t,w;\xb). 
$$  
\end{cor} 

In the case of $n =2m$ we can use \eqref{eq:projEven} instead of 
\eqref{eq:proj} in the last step of the proof of Proposition \ref{prop:4.4}
to get an expression for $\proj_{l,3}^\mu f$. The corresponding expression
for the partial sum is the following result: 

\begin{prop} \label{prop:4.4Even} 
For $m \ge 0$, 
$$
S_{2m,3}^\mu f(\xb) = \frac{1}{2m+1}\sum_{\nu =0}^{2m} \int_{-1}^1\int_{-1}^1
 \CR_{\phi_\nu}^\mu (f;t,w) \Phi_{2m}(\phi_\nu,t, w; \xb) dt \sqrt{1-w^2} dw.
$$
\end{prop}

From such an expression of $S_{n,3}^\mu$ we naturally want to derive an 
algorithm as in the 2D case. However, there is a problem when we use 
quadrature formula. Indeed, in order to obtain an algorithm, we need to 
discretize the integrals
\begin{equation}\label{eq:integ}
  \int_{-1}^1\int_{-1}^1
 \CR_{\xi_\nu}^\mu(f;t,w)\Phi_{n}^\mu (\xi_\nu,t,w;\xb)dt\sqrt{1-w^2}dw 
\end{equation}
in $S_{n,3}^\mu f$ by a quadrature formula. We can use, for example, the
quadrature \eqref{eq:gaussian} of precision $2n$, which we denote by
$$
   \int_{-1}^1 f(t) (1-t^2)^\alpha dt \approx 
       \sum_{k=0}^{n} \lambda_{k,n}^\alpha f(t_{k,n}^\alpha)
$$
to emphasis the dependence of $t_{k,n}$ and $\lambda_{k,n}$ on the weight
function. If we follow the 2D case, then the equation \eqref{eq:RQ3} 
indicates that we should apply the quadrature with respect to $(1-t^2)^\mu$ 
in $t$ variable, and apply the quadrature with respect to $(1-w^2)^{\mu+1/2}$ 
in $w$ variable. The result of using these quadrature formulas gives the 
following: 

\begin{algo} \label{algo:B3}
Let $\mu \ge 0$. For $n \ge 0$, $\xb = (x_1,x_2,x_3) \in B^3$, 
$$
\CB_n^\mu (f; \xb) = \sum_{\nu = 0}^{n} \sum_{j=0}^{n} \sum_{k=0}^{n} 
   \CR_{\phi_\nu}^\mu (f; t_{j,n}^\mu, t_{k,n}^{\mu+1/2})T_{j,k,\nu}(\xb) 
$$
where 
$$
 T_{j,k,\nu}^\mu (\xb) = \frac{\lambda_j^\mu \lambda_k^{\mu+1/2}}{n+1} 
        \Phi_{n}^\mu (\xi_\nu, t_{j,n}^\mu, t_{k,n}^{\mu+1/2}; \xb). 
$$
\end{algo} 
 
However, this is likely not an accurate algorithm. The problem is that the
operator $\CB_n^\mu$ does not preserve polynomials of degree $n$. In fact, 
in order that $\CB_n P = P$ for $P \in \Pi_n^3$, we need the discretization
of the integrals \eqref{eq:integ} to be exact whenever $f$ is a polynomial
of degree at most $n$. The function
$$
F_\mu(t,w):=(1-t^2)^{-\mu}(1-w^2)^{-\mu}
      \CR_{\xi_\nu}^\mu (f;t,w) \Phi_n^\mu(\xi_\nu, t,w;\xb)
$$
is a polynomial of degree $2n$ in variable $t$ whenever $f$ is a polynomial
of degree $n$ by the definition of $\Phi_n^\mu$ and Proposition \ref{prop:4.3},
so that the discretization in $t$ variable is exact. However, the function 
$F_\mu(t,w)$ is not a polynomial in $w$ variable. By the definition of 
$Q_{l,k,j}$ in \eqref{eq:Qlkj}, the equation \eqref{eq:RQ3} shows that 
$F_\mu(t,w)$ with $f = Q_{l,k,j}$ contains $(1-w)^{k/2}\wt 
C_{l-k}^{k+\mu+1}(w)$, which is not a polynomial in $w$ variable if $k$ is
odd. The formula of $\Phi_n^\mu(\xi_\nu, t,w; \xb)$ shows that it is a 
sum of functions, which is also not a polynomial. This means that the 
quadrature will not be exact and polynomials are not preserved by $\CB_n^\mu$. 

An algorithm should have high convergence order if it preserves polynomials 
up to certain degrees. The fact that $B_n^\mu f$ does not preserves polynomials
means that the convergence of the algorithm may not be as desirable. 

%For a given function $f$, the approximation $\CB_{2m}f$ uses the set of 
%the attenuated Radon projections
%$$
%  \left\{ \CR_{\phi_\nu}^\mu (f; t_{j,n}^\mu,t_{k,n}^{\mu+1/2}):
%         0 \le \nu \le 2m, 0 \le j \le 2m, 0 \le k \le 2m \right\}
%$$
%of $f$, which consists of projections on $2m+1$ parallel planes defined by
%$x_3 = t_{k,n}^{\mu+1/2}$. The intersection of such an plane with $B^3$ is
%a disk, call it $B_3^k$. There are $2m_1$ equally spaced directions along the 
%circumference of the disk $B_k^3$ (specified by $\xi_\nu$) and there are 
%$2m+1$ parallel lines (specified by $t_{k,n}^\mu$) in each direction. 

%In the case of $\mu =1/2$, we are dealing with the usual Radon projections
%with weight function $W_{1/2}(x) =1$. The algorithm is easy to implement
%since $T_{j,k,\mu}(\xb)$ are the fixed polynomials and can be saved on the 
%hard drive before hand. The operator $\CB_{2m} f$ also satisfies the following
%property: 

%%%%%%%%%%%%%%%%%%%%%%%
%%%%%% Section 5  %%%%%
%%%%%%%%%%%%%%%%%%%%%%%

\section{Reconstruction and Approximation on the cylinder domain}
\setcounter{equation}{0}

In contrast to the unit ball in $\RR^3$, the reconstruction algorithm on
a cylinder domain works well. Let $L > 0$ and let $B_L$ be the cylinder domain
defined by
$$
   B_L = B^2 \times [0, L] =\{(x,y,z): (x,y) \in B^2, 0 \le z \le L\}.
$$  
We will show that the partial sum operator of the orthogonal expansions on 
$W_L$ admits an expression that relates to Radon data and use it to get a
reconstruction algorithm. 

Let $W_\mu$ be defined as in \eqref{eq:W-mu}. Let $W_{\mu,L}$ be the weight 
function 
$$
  W_{\mu,L}(x,y,z) = W_\mu(x,y) W_L(z), \qquad (x,y,z) \in B_L.
$$
We retain the notation  $\CR_\phi^\mu(g;t)$ for the attenuated Radon 
projection of a function $g: B^2 \mapsto \RR$, as defined in 
\eqref{eq:R-mu}. For a fixed $z$ in $[0,L]$, we define
\begin{equation}\label{eq:Radon3DL}
 \CR_\phi^\mu(f(\cdot,\cdot,z); t) := \int_{I(\phi,t)} f(x,y,z) 
    W_\mu(x,y)dx dy,
\end{equation}
which is the attenuated Radon projection of $f$ in a disk that is 
perpendicular to the $z$-axis.  

We consider the orthogonal polynomials with respect to the inner product 
\begin{equation} \label{eq:inner}
 \langle f, g\rangle_{B_L}  =  \frac{1}{\pi} \int_{B_L} f(x,y,z) g(x,y,z)
       W_{\mu,L}(x,y,z)\,dx\, dy \,dz.
\end{equation}
Let $\CV_n^3(W_{\mu,L})$ denote the subspace of orthogonal polynomials of 
degree $n$ on $B_L$ with respect to the inner product \eqref{eq:inner};  
that is, $P \in \CV_n^3(W_{\mu,L})$ if $\langle P, Q\rangle_{B_L} =0$ for 
all polynomial $Q \in \Pi_{n-1}^3$. 

Let $p_k$ be the orthonormal polynomial of degree $n$ with respect to 
$W_L$ on $[0,L]$ and let $\{P_j^k(x,y): 0 \le j \le k\}$ denote an orthonormal
basis of $\CV_k^2(W_\mu)$. Since $W_{\mu,L}$ is a product on a product domain,
the following proposition is obvious. 

\begin{prop} \label{prop:5.1}
An orthonormal basis for $\CV_l^3(W_{\mu,L})$ is given by 
$$
\PP_l=\left\{P_{l,k,j}^\mu: 0 \le j \le k \le n\right\}, \qquad 
    P_{n,k,j}^\mu(x,y,z) = P_j^k(x,y)p_{n-k}(z).
$$
In particular, the set $\{\PP_l: 0 \le l \le n \}$ is an orthonormal basis 
for $\Pi_n^3$. 
\end{prop}

For $f \in L^2(W_{\mu,L};B_L)$, the Fourier coefficients of $f$ with 
respect to the orthonormal system $\{\PP_l: l \ge 0\}$ are given by 
$$
 \wh f_{l,k,j}^\mu = a_\mu \int_{B_L} f(\xb)P_{l,k,j}^\mu(\xb) 
      W_{\mu,L}(\xb)d\xb,  \quad 0 \le j \le k \le  l.
$$ 
Let $S_{n,L}^\mu f$ denote the Fourier partial sum operator, 
$$
S_{n,L}^\mu f(\xb) = \sum_{l=0}^n \sum_{k=0}^l \sum_{j=0}^k 
   \wh f_{l,k,j}^\mu P_{l,k,j}^\mu(\xb). 
$$
Just like its counterpart in two variables,  this is a projection operator.
The following is an analogue of Theorem \ref{thm:partial-sum} for the
cylinder domain $B_L$. 

\begin{thm} \label{thm:partial-sum3D}
For $n \ge 0$, 
\begin{align} \label{eq:Sn3D}
  S_{n,L}^\mu f(\xb)  = \frac{1}{n+1} \sum_{\nu=0}^{n} 
    a_\mu \int_{-1}^1 \int_0^L \CR_{\xi_\nu}^\mu (f(\cdot,\cdot,w);t)
      \Phi_n^\mu(\xi_\nu, w,t;\xb)  W_L(w)dw\,dt 
\end{align}
where
\begin{equation} \label{eq:Phi3D}
  \Phi_n^\mu (\xi, w,t;\xb) = \sum_{k=0}^{n} \frac{k+\mu+1/2}{\mu+1/2} 
     D_k^{\mu+1/2}(\xi, t;x_1,x_2) \sum_{l=0}^{n-k} p_l(w)p_l(x_3).
\end{equation}
\end{thm} 

\begin{proof}
By the definition of $\wh f_{l,k,j}^\mu$ we can write 
$$
 \wh f_{l,k,j}^\mu = a_\mu \int_{B^2} f_{l-k} (x,y) P_j^k(x,y)W_\mu(x,y)dxdy 
$$
where
$$
   f_{l-k}(x,y) : = \int_0^L f(x,y,w) p_{l-k}(w) W_L(w) dw, 
     \qquad  l \ge k \ge 0. 
$$
Consequently, by the definition of $\proj_k^\mu$ in \eqref{eq:Pintegral},
it follows that 
$$
 S_{n,L}^\mu f(\xb) = \sum_{l=0}^n\sum_{k =0}^l \proj_k^\mu (f_{l-k}; x_1,x_2)
    p_{l-k}(x_3).
$$
We can then use the expression \eqref{eq:proj} for $\proj_k^\mu f$ and the 
fact that 
$$
  \CR_{\xi}^\mu (f_{l-k};t) = \int_0^L \CR_\xi^\mu (f(\cdot,\cdot,w); t)
      p_{l-k}(w) W_L(w)dw
$$
to complete the proof. 
\end{proof}

In the case of $n =2m$, we can use \eqref{eq:projEven} in place of 
\eqref{eq:proj} in the proof. The result is the following proposition which
has appeared in \cite{X05} when $\mu =1/2$.

\begin{prop}
For $m \ge 0$, 
\begin{align} \label{eq:S2m3D}
 & S_{2m,L}^\mu f(\xb) \\
 & \quad  = \frac{1}{2m+1} \sum_{\nu=0}^{2m} 
    a_\mu \int_{-1}^1 \int_0^L \CR_{\xi_\nu}^\mu (f(\cdot,\cdot,w);t)
      \Phi_{2m}^\mu(\phi_\nu, w,t;\xb)  W_L(w)dw\,dt. \notag
\end{align}
\end{prop} 

From the expression \eqref{eq:Sn3D} or \eqref{eq:S2m3D} of $S_{n,L}^\mu f$,
we can apply a quadrature formula to get a reconstruction algorithm on $B_l$
for the attenuated Radon data. In \cite{X05} the weight function
$W_L$ is chosen to be the Chebyshev weight function
$$
   W_L(z) = \frac{1}{\pi} \frac{1}{\sqrt{z(L-z)}}, \qquad z \in [0, L],
$$
normalized to have integral $1$ on $[0,L]$. The reason for this choice is
that the Gaussian quadrature formula takes a simple form
\begin{equation}\label{eq:quadratureT}
 \int_{0}^L g(z) W_L(z) dz \approx \frac{1}{n+1}\sum_{j=0}^n g(z_i), \qquad
    z_i = \frac{1}{2} \left(1+ \cos \tfrac{2j+1}{2n+2}\right), 
\end{equation}
which is of precision $2n+1$. We can apply this quadrature for the integral
with respect to $w$ and use the quadrature \eqref{eq:gaussian} for the 
integral with respect to $t$ in \eqref{eq:Sn3D} or \eqref{eq:S2m3D}. The
result is the following algorithm:

\begin{algo} \label{algo:3DBL}
Let $\mu \ge 0$ and let $\gamma_{\mu,j,i}=\CR_{\xi_\nu}^\mu(f(\cdot,\cdot,z_i);
 t_{j,n})$. For $n \ge 0$
\begin{equation}
 \CB_{n,L}^\mu (f;\xb) = \sum_{\nu = 0}^n \sum_{j= 0}^n \sum_{i = 0}^n 
    \gamma_{\nu,j,i} T_{\nu,j,i}(\xb)
\end{equation}   
where 
$$
 T_{\nu,j,i}(\xb) = \frac{a_\mu \lambda_{j,n}}{n+1}
  (1-t_{j,n}^2)^{-\mu} \Phi_n^\mu(\xi_\nu, z_i, t_{j,n};\xb).
$$
\end{algo} 

Like the algorithms in the previous sections, this algorithm produces a 
polynomial as an approximation to the function. It does preserve 
polynomials of lower degrees.

\begin{thm}
The operator $\CB_{n,L}^\mu$ is a projection operator on $\Pi_n^3$. In
other words, $\CB_n f \in \Pi_n^3$ and $\CB_{n,L}(f) = f$ if $f \in \Pi_n^3$. 
\end{thm}  

\begin{proof}
Let $P_{n,k,j}^\mu$ be defined as in Proposition \ref{prop:5.1}. It follows
from the definition in \eqref{eq:Radon3DL} that 
$
\CR_\phi^\mu(P_{l,k,j}^\mu(\cdot,\cdot,w);t)=\CR^\mu_\phi(P_j^k;t)p_{l-k}(w). 
$ 
Consequently, it follows from \eqref{eq:Marr} that %for $P\in \CV_l^2$,
%$
%\CR_\phi^\mu(P(\cdot,\cdot,z);t)= b_\mu(1-t^2)^\mu
%         \frac{C_k^{\mu+1/2}(t)}C_k^{\mu+1/2}(1)} p_{l-k}(z). 
%$$
$\CR_\phi^\mu(P(\cdot,\cdot,w);t)/(1-t^2)^{\mu}$ is a polynomial of degree
$n$ in both $t$ variable and $w$ variable whenever $P \in \Pi_n^3$. By its
definition in \eqref{eq:Phi3D}, the function $\Phi^\mu(\xi, w, t;\xb)$ is 
evidently a polynomial of degree $n$ in both $t$ and $w$ variables. Hence,
we can apply \eqref{eq:quadratureT} for $w$ variable and apply the quadrature 
\eqref{eq:gaussian} of precision $2n$ to $t$ variable, which are exact on
$(1-t^2)^{-\mu}\CR_\phi^\mu(P(\cdot,\cdot,w);t)\Phi^\mu(\xi, w, t;\cdot)$.
\end{proof} 

The approximation process in Algorithm \ref{algo:3DBL} uses the attenuated
Radon data 
$$
  \left\{ \CR_{\xi_\nu}^\mu (f(\cdot,\cdot,z_i);t_{j,n}): 
     0 \le \nu \le n, \,\, 0 \le j \le n,\,\, 0 \le i \le n \right\},
$$
which consists of Radon projections on $n+1$ disks that are parallel to 
the $z$-axis. In other words, it consists of reconstructions of the function
on $n+1$ planes. 

In the case of $n = 2m$ and $\mu$ is an half integer, we can also use the
quadrature \eqref{eq:gaussianU} to derive a more explicit algorithm as in 
Algorithm \ref{algo:N/2}. Such an algorithm is given in \cite{X05} for 
$\mu =1/2$. We shall not elaborate further.  

%\medskip

%{\it Acknowledgment}: The author thanks Dr. Christoph Hoeschen and 
%Dr. Oleg Tischenko for stimulating discussions about the computer tomography. 


\begin{thebibliography}{99}
 
\bibitem{BO} 
        T. Bortfeld and U. Oelfke, 
        Fast and exact 2D image reconstruction by means of Chebyshev
        decomposition and backprojection,
        \textit{Phys. Med. Biol.}. \textbf{44} (1999), 1105-1120.

\bibitem{DX}
        C. F. Dunkl and Yuan Xu,
        \textit{Orthogonal polynomials of several variables},
        Cambridge Univ. Press, 2001. 

\bibitem{F}
        D. Finch, 
        The attenuated x-ray transform: recent developments,
        in \textit{Inside out: inverse problems and applications}, 
        47--66, Math. Sci. Res. Inst. Publ., \textbf{47}, 
        Cambridge Univ. Press, Cambridge, 2003.

\bibitem{KS}
        A. C. Kak and M. Slaney,
        \textit{Principles of Computerized Tomographic Imaging},
        IEEE Press, New York, 1988; Reprint as Classics in Applied Mathematics,
        \textbf{33}.  SIAM, Philadelphia, PA, 2001.

\bibitem{LS} 
        B. Logan and L. Shepp, 
        Optimal reconstruction of a function from its projections, 
        \textit{Duke Math. J.} \textbf{42} (1975), 645-659.    
       
\bibitem{Marr} 
        R. Marr, 
        On the reconstruction of a function on a circular domain from 
        a sampling of its line integrals, 
        \textit{J. Math. Anal. Appl.}, \textbf{45} (1974), 357-374.
 
\bibitem{N} 
        F. Natterer,
        \textit{The mathematics of computerized tomography},
        Reprint of the 1986 original. Classics in Applied Mathematics, 
        32. SIAM, Philadelphia, PA, 2001. 

\bibitem{N2} 
        F. Natterer, 
        Inversion of the attenuated Radon transform,
        \textit{Inverse Problems}, \textbf{17}  (2001),  no. 1, 113--119.

\bibitem{NW} 
        F. Natterer and F. W\"ubbeling, 
        \textit{Mathematical Methods in Image Reconstruction},
        SIAM, Philadelphia, PA, 2001. 
 
\bibitem{No}
        R. G. Novikov,
        An inversion formula for the attenuated X-ray transformation,
        \textit{Ark. Mat.} \textbf{40} (2002), 145--167.

\bibitem{Sz} 
	G. Szeg\H{o},
	\textit{Orthogonal Polynomials},  
	Amer. Math. Soc. Colloq. Publ. Vol.23, Providence, 4th edition,
        1975.
  
\bibitem{X98}
        Yuan Xu,
         Orthogonal polynomials and cubature formulae on spheres and on
         balls, 
         \textit{SIAM J. Math. Anal.} \textbf{29} (1998), 779--793.
 
\bibitem{X99}
        Yuan Xu,
        Summability of Fourier orthogonal series for Jacobi weight on a 
        ball in $\RR^d$, \textit{Trans. Amer. Math. Soc.}  \textbf{351} (1999),
        2439-2458.

\bibitem{X00}
        Yuan Xu,
        Funk-Hecke formula for orthogonal polynomials on spheres and
        on balls, 
        \textit{Bull. London Math. Soc.} \textbf{32} (2000), 447-457.
 
\bibitem{X01}
        Yuan Xu,
        Representation of reproducing kernels and the Lebesgue constants
        on the ball,       
        \textit{J. Approx. Theory} \textbf{112} (2001), 295-310.
 
\bibitem{X05}
        Yuan Xu,
        A new approach to the reconstruction of images from Radon projections,
        \textit{Adv. in Applied Math.}, accepted for publication. 
        %submitted, Feb. 2005.
        
\bibitem{XTC}
        Yuan Xu, O. Tischenko, and C. Hoeschen, 
        New tomographic reconstruction algorithms,
        \textit{submitted}, 2005.

\bibitem{XTC2}
        Yuan Xu, O. Tischenko, and C. Hoeschen, 
        A new reconstruction algorithm for Radon Data, 
        \textit{SPIE Proceedings of Medical Imaging, 2006},
        to appear.
   
\bibitem{Z}
         A. Zygmund, 
         \textit{Trigonometric Series}, 
         Cambridge Univ. Press, 1959. 
            
\end{thebibliography}
\end{document}